\documentclass{article}
\usepackage{amsfonts}
\usepackage{amsmath}
\usepackage{amssymb}
\usepackage{amsthm}
\usepackage{fullpage}
\usepackage{cite}

\newtheorem{theorem}{Theorem}

\newtheorem{corollary}[theorem]{Corollary}

\newtheorem{lemma}[theorem]{Lemma}
\newtheorem{proposition}[theorem]{Proposition}

 \DeclareMathOperator{\supp}{supp}
\DeclareMathOperator{\curl}{curl} \DeclareMathOperator{\diver}{div}
\numberwithin{equation}{section} \numberwithin{theorem}{section}
\allowdisplaybreaks[1]

\begin{document}

\title{On the convergence rate of the Euler-$\alpha$, an inviscid
second-grade complex fluid, model to the Euler equations}
\author{Jasmine S.~Linshiz$^{1a}$ and Edriss S.~Titi$^{1,2}$}
\date{November 8, 2009}

\maketitle

\begin{center}
$^{1}$\textit{Department of Computer Science and Applied Mathematics \\[0pt]
Weizmann Institute of Science \\[0pt]
Rehovot 76100, Israel}\\[0pt]
$^{a}$jasmine.tal@weizmann.ac.il \\[0pt]
$^{2}$\textit{Department of Mathematics and\\[0pt]
Department of Mechanical and Aerospace Engineering \\[0pt]
University of California \\[0pt]
Irvine, CA 92697-3875, USA} \\[0pt]
etiti@math.uci.edu \textit{and} edriss.titi@weizmann.ac.il
\end{center}

\noindent \textbf{Keywords:} inviscid regularization of Euler
equations; Euler-$\alpha$; second-grade non-Newtonian fluid; vortex
patch.

\noindent \textbf{Mathematics Subject Classification:} 76B03, 35Q35, 76B47.

\begin{abstract}
We study the convergence rate of the solutions of the incompressible Euler-$%
\alpha $, an inviscid second-grade complex fluid, equations to the
corresponding solutions of the Euler equations, as the
regularization
parameter $\alpha $ approaches zero. First we show the convergence in $H^{s}$%
, $s>n/2+1$, in the whole space, and that the smooth Euler-$\alpha $
solutions exist at least as long as the corresponding solution of
the Euler equations. Next we estimate the convergence rate for
two-dimensional vortex patch with smooth boundaries.
\end{abstract}

\section{\label{sec:intro}Introduction}

The equations of motion for a visco-elastic second-grade non-Newtonian
complex fluid are given by the system (see, e.g., \cite{a_DF74,a_DR95,Truesdell})%
\begin{align}
& \frac{\partial v}{\partial t}-\nu \Delta {u}+\left( {u}\cdot \nabla
\right) {v}+{v}_{j}\nabla {u}_{j}+\nabla p=0,  \label{grp:vis2gradeEq} \\
& v=\left( 1-\alpha ^{2}\Delta \right) {u},  \notag \\
& \nabla \cdot u=0,  \notag \\
& v(x,0)=v^{in}(x),  \notag
\end{align}%
where the fluid velocity field, $v$, and the pressure, $p$, are the
unknowns; $v^{in}$ is the given initial velocity field, $\nu $ is kinematic
viscosity, $\alpha >0$ is a material parameter which represents the elastic
response of the fluid, and we use Einstein's summation convention.The
inviscid version of this model, i.e., when $\nu =0$, is mathematically
identical to the Euler-$\alpha $ [also known as the Lagrangian-averaged
Euler-$\alpha $] model, which was independently introduced and derived in
the Euler-Poincar\'{e} variational framework in \cite{a_HMR98a,a_HMR98b}. In
this variational theory the parameter $\alpha $ is interpreted as a spatial
filtering scale of the velocity field $v$.

In \cite{a_CFHOTW98,a_CFHOTW99, a_CFHOTW99_ChanPipe,a_FHT01,a_FHT02} the
corresponding Navier-Stokes-$\alpha $ (NS-$\alpha $) [also known as the
viscous Camassa-Holm equations or the Lagrangian-averaged Navier-Stokes-$%
\alpha $ (LANS-$\alpha $)] model, which is a regularization of the
Navier-Stokes equations (NSE), was obtained by introducing an
appropriate {\it ad hoc} viscous term into the Euler-$\alpha $
equations, that is, by adding the
viscous term $-\nu \Delta {v}$, instead of $-\nu \Delta {u}$ in %
\eqref{grp:vis2gradeEq}. While the question of global regularity of
three-dimensional (3D) viscoelastic model \eqref{grp:vis2gradeEq} is still a
challenging open problem, the 3D NS-$\alpha $ model is proven to be globally
well posed \cite{a_FHT02}. The extensive research of the $\alpha $-models
(see, e.g., \cite%
{a_CT09,a_BLT09,a_BLT08,a_JNTX09,a_GKT08,a_CHMZ99,a_HT05,a_CHOT05,a_ILT05,a_VTC05,a_CTV07,a_FHT02,a_FHT01,a_CFHOTW99,a_CFHOTW99_ChanPipe, a_CFHOTW98,a_MKSM03,a_LL06,a_LL03,b_BIL06,a_L06,a_GKT08,a_GH03,a_GH06,a_CLT06,a_BFR80,a_HN03,a_CHT04,a_CFR79,a_LT07,a_GHMP08}%
) stems, on the one hand, from the close agreement of their steady state
solutions to averaged empirical data, for a large range of huge Reynolds
numbers, for turbulent flows in infinite channels and pipes \cite%
{a_CFHOTW99,a_CFHOTW99_ChanPipe, a_CFHOTW98}. On the other hand, the $\alpha
$-models, for small values of the parameter $\alpha $, can also be viewed as
numerical regularizations of the original, Euler or Navier-Stokes, systems.
The main practical question arising is that of the applicability of these
regularizations to the correct predictions of the underlying flow phenomena.
In particular, it becomes important to investigate the problem of
convergence of the $\alpha $-models, as the regularization parameter $\alpha
$ approaches zero. This problem has been studied in various contexts. In
\cite{a_FHT02} the convergence, as $\alpha \rightarrow 0$, of a subsequence
of the weak solutions of the three dimensional (3D) NS-$\alpha $ equations
to a Leray-Hopf weak solution of the 3D NSE equations is shown, for the case
of periodic boundary conditions. Similar results are also reported in \cite%
{a_LT07} concerning the MHD-$\alpha $ model. In \cite{a_CTV07,a_VTC07} it is
shown that the trajectory attractors of the 3D Leray-$\alpha $ and NS-$%
\alpha $ models converge to the trajectory attractor of weak solutions of
the NSE as $\alpha \rightarrow 0$. The authors of \cite{a_CGKTW08} obtained
a rate of convergence of the solutions of the 3D NS-$\alpha $ equations with
periodic boundary to the solutions of the Navier-Stokes equations, as $%
\alpha \rightarrow 0$, for small initial data in Besov-type function spaces,
for which global existence and uniqueness of solutions can be established.
Recently, the convergence rates of solutions of various two-dimensional (2D)
$\alpha $-regularization models, subject to periodic boundary conditions,
toward solutions of the exact Navier-Stokes equations, as $\alpha
\rightarrow 0$, have been studied in \cite{a_CT09}. For 2D Euler-$\alpha $
regularization the authors of \cite{a_BLT09} show the convergence of a
subsequence of the weak solutions of the Euler-$\alpha $ equations with a
distinguished sign vortex sheet initial data (a Radon measure supported on a
curve) to a solution of the 2D Euler equations, as $\alpha \rightarrow 0$.
We elaborated on this result below. It is also worth mentioning that the
problem of weak convergence of solutions of the viscous second-grade
equations \eqref{grp:vis2gradeEq} with $L^2$ initial data to a solution of
the Navier-Stokes equations, as $\alpha \to 0$, is treated in \cite{a_I02a,
a_I02}.

In this paper we follow the above mentioned philosophy proposed in \cite%
{a_LT07} and consider the Euler-$\alpha $ model as a numerical
inviscid regularization of the Euler equations. We study the
convergence rate of the solutions of the 2D and 3D Euler-$\alpha $
equations to the corresponding solutions of the Euler equations for
smooth initial data in the whole space, as the regularization
parameter $\alpha $ approaches zero. In the 2D case we also
investigate the convergence rate of the solutions of the
Euler-$\alpha $ to the corresponding solutions of the Euler
equations for vortex patch initial data. This program is analogues
to the established results concerning the rate of convergence of the
NSE to the Euler equations, as the viscosity $\nu \to 0$, see
\cite{a_K72,a_M07,a_AD04,a_CW95}. However, in the 2D case there is
an advantage of the Euler-$\alpha $ regularization over the
Navier-Stokes equations regularization since the former regularizes
the solution by transporting the vorticity with a smooth vector
field, and hence it preserves the structure of the vortex patch and
vortex sheet while regularizing the motion.

The incompressible Euler equations are given by

\begin{align}
& \frac{\partial v}{\partial t}+(v\cdot \nabla )v+\nabla p=0,
\label{grp:EulerEq} \\
& \nabla \cdot v=0,  \notag \\
& v(x,0)=v^{in}(x),  \notag
\end{align}%
where $v$, the fluid velocity field, and $p$, the pressure, are the
unknowns, and $v^{in}$ is the given initial velocity field. For results
concerning the Euler equations, see \cite{b_MP94,b_MB02}, and for recent
surveys, see \cite{a_BT07,a_C07}.

In $\mathbb{R}^{2}$, the vorticity formulation of Euler equations is
obtained by taking a curl of \eqref{grp:EulerEq} and is given by
\begin{align}
& \frac{\partial q}{\partial t}+\left( v\cdot \nabla \right) q=0,
\label{grp:EulerEqVortForm} \\
& v=K\ast q,  \notag \\
& q(x,0)=q^{in}(x),  \notag
\end{align}%
where $K\left( x\right) =\frac{1}{2\pi }\nabla ^{\perp }\log \left\vert
x\right\vert $, \mbox{$q=\curl v$} is the vorticity, and $q^{in}$ is the
given initial vorticity.

Yudovich \cite{a_Y63} obtained the existence and uniqueness of weak
solutions of the 2D incompressible Euler equations for initially bounded
vorticity (see, also, \cite{a_B72,a_K92} for an alternative proof, \cite%
{a_V99} for an improvement with vorticity in a class slightly larger than $%
L^{\infty }$, and \cite{a_T04} for review of relevant two-dimensional
results). In particular, the problem of evolution of vortex patches, where
the vorticity is a multiple of the characteristic function of a bounded
domain, has a unique global solution, and it was proved in \cite{a_C93} (see
also \cite{a_BC93}) that $C^{1,\gamma }$, $\gamma >0$, boundaries of the
patches remain $C^{1,\gamma }$ for all times. In \cite{a_CW95} it was shown
that the $L^{2}$ norm of the difference between the solutions of NSE and the
corresponding solution of the Euler system for such initial data converges
to zero, as the kinematic viscosity $\nu \rightarrow 0$; even though none of
the solutions is in $L^{2}$. The $\left( \nu t\right) ^{1/2}$ rate of
convergence of \cite{a_CW95} was improved to $\left( \nu t\right) ^{3/4}$ in
\cite{a_AD04}, due to the fact that the vorticity of the vortex patch with $%
C^{1,\gamma }$ boundary is in fact in a Besov space $\dot{B}_{2,\infty
}^{1/2}$, see also \cite{a_M07} for a simpler proof and an extension to $%
\mathbb{R}^{3}$.

As we have mentioned above the Euler-$\alpha $ model \cite%
{a_CFHOTW99_ChanPipe,a_HMR98a,a_HMR98b,a_H02_pA,a_MS03,a_C01} is an inviscid
regularization of the Euler equations \eqref{grp:EulerEq}, which is given by
the system
\begin{align}
& \frac{\partial v^{\alpha }}{\partial t}+\left( {u}^{\alpha }\cdot \nabla
\right) {v}^{\alpha }+{v}_{j}^{\alpha }\nabla {u}_{j}^{\alpha }+\nabla
p^{\alpha }=0,  \label{grp:EulerAlphaEq} \\
& v^{\alpha }=\left( 1-\alpha ^{2}\Delta \right) u^{\alpha },  \notag \\
& \nabla \cdot u^{\alpha }=\nabla \cdot v^{\alpha }=0,  \notag \\
& v^{\alpha }(x,0)=v^{in,\alpha }(x),  \notag
\end{align}%
Here $u^{\alpha }$ represents the \textquotedblleft filtered" fluid velocity
vector, $p^{\alpha }$ is the \textquotedblleft filtered" pressure, $\alpha
>0 $ is a regularization length scale parameter representing the width of
the filter. Observe that for $\alpha =0$ one recovers, formally, the Euler
equations \eqref{grp:EulerEq}. The vorticity of the 2D Euler-$\alpha $ model %
\mbox{$q^{\alpha}=\curl v^{\alpha}$} obeys the equations
\begin{align}
& \frac{\partial q^{\alpha }}{\partial t}+\left( u^{\alpha }\cdot \nabla
\right) q^{\alpha }=0,  \label{grp:EulerAlphaEqVortForm} \\
& u^{\alpha }=K^{\alpha }\ast q^{\alpha },  \notag \\
& q^{\alpha }(x,0)=q^{in,\alpha }(x).  \notag
\end{align}%
The smoothed kernel is $K^{\alpha }=G^{\alpha }\ast K$, where $G^{\alpha
}\left( x\right) =\frac{1}{\alpha ^{2}}\frac{1}{2\pi }K_{0}\left( \frac{%
\left\vert x\right\vert }{\alpha }\right) $ is the Green function associated
with the Helmholtz operator $\left( I-\alpha ^{2}\Delta \right) $ (see,
e.g., \cite{b_P02}), the function $K_{0}$ is a modified Bessel function of
the second kind of order zero (see, e.g.,\cite{b_W44}).

The 2D Euler-$\alpha $ equations were studied in \cite{a_OS01}, where it has
been shown that there exists a unique global weak solution to the Euler-$%
\alpha $ equations with initial vorticity in the space of Radon measures on $%
{\mathbb{R}}^{2}$, with a unique Lagrangian flow map describing the
evolution of particles. We remark, however, that the question of global
existence of weak solutions for the three-dimensional (3D) Euler-$\alpha $
equations is still an open problem. In \cite{a_JNTX09} the global existence
of weak solutions is shown for the 3D axisymmetric Euler-$\alpha $ equations
without swirl, for initial vorticity being a finite Radon measure with
compact support, also, the global existence and uniqueness is established
for the compactly supported vorticity in $L^{p}$, $p>3/2$, see also \cite%
{a_BR04}.

As in the 2D Euler case, the vortex patch is transported by the flow under
evolution of Euler-$\alpha $ equations. The following result is the essence
of Theorem \ref{thm:GlobalEx} stated below. If the initial vorticity is a
multiple of the characteristic function of a simply connected bounded domain
$\Omega ^{in}$, with a certain technical condition on the boundary being a
simple curve, and the boundary $\partial \Omega ^{in}$ is in either one of
the following spaces: $\mathrm{Lip}$, or $C^{1,\beta }$, $0\leq \beta \leq 1$%
, or $C^{2,\beta }$, $0\leq \beta <1$, or $C^{n,\beta }$, $n\geq 3$, $%
0<\beta <1$, then the boundary of the vortex patch remains in the
corresponding space for all time.

The outline of the paper is as follows. In Section \ref%
{sec:classical_solutions} we study the convergence rate of the solutions of
the Euler-$\alpha $ equations to the solution of the Euler equations for
strong solutions that belong to the Sobolev space $H^{m}\left( \mathbb{R}%
^{n}\right) $, $m>n/2+1$, for $n=2,3$. We show that interval of existence of
Euler-$\alpha $ solutions contains the interval of existence of the
corresponding Euler solution, and that the convergence is uniform in time,
for time intervals compactly contained in $\left[ 0,T^{\ast }\right) $,
where $T^{\ast }$ is the time of existence of the solution of the Euler
equations. We also show that in the $H^{m-2}$ norm the solution of the Euler-%
$\alpha $ equations differs from the solution of the Euler equations
by order $\alpha ^{2}$. It is worth mentioning that this result
corresponds to the $\nu t$ convergence rate of the solutions of the
NSE to the one of the Euler equations, for the inviscid limit of the
classical solutions of the NSE equations in the whole space see
\cite{a_K72, a_K75,a_M07,a_M67,a_S71}. The issue of inviscid limit
of the NSE in domains with physical boundaries, subject to the
no-slip Dirichlet boundary conditions, is a very important open
problem, for both theoretical study and applications. The problem
emerges first from the boundary layer, which appears because we can
not impose a Dirichlet boundary condition for the Euler equation,
then the nonlinear advection term of the Navier--Stokes equations
may propagate this instability inside the domain. Very few
mathematical results are available for this very unstable situation.
One of the most striking results in this direction is a theorem of
Kato \cite{a_K84}, see also \cite{a_TW97,a_W01}.

In Section \ref{sec:vortex_patch} we study the convergence for the vortex
patch problem. Specifically, we show the convergence, as $\alpha \rightarrow
0$, of the $L^{2}$ norm of the difference between the solutions of the Euler-%
$\alpha $ equations and the solution of the Euler equations for the vortex
patch initial data with vorticity being a characteristic function of a
simply connected bounded domain with $C^{1,\gamma }$, $\gamma >0$, boundary,
even though neither of the solutions are in $L^{2}$. The convergence rate is
of order $\left( \alpha ^{2}\right) ^{3/4}$, which corresponds to the
optimal convergence rate of the difference between the solutions of the NSE
and the Euler equations which is of order $\left( \nu t\right) ^{3/4}$ \cite%
{a_AD04}.

We remark that, ideally, we would like to compute the rate of convergence in
dimensionless units, however due to the absence of typical length scale in $%
\mathbb{R}^{n}$, the above rates of convergence of the NSE to the Euler
equations are given as powers of $\left( \nu t\right) $, which has the units
of length square, see, e.g., \cite{a_K72,a_M07,a_AD04,a_CW95}. Similarly, in
our case, the rate of convergence involves the parameter $\alpha ^{2}$.
Observe, that one can artificially cook up a length scale from the initial
value , e.g., $\left\Vert v^{in}\right\Vert _{L^{2}}\left\Vert
q^{in}\right\Vert _{L^{2 }}^{-1} $ in the case of $H^{m}\left( \mathbb{R}%
^{n}\right) $, $m>n/2+1$, solutions, or $\left\Vert q^{in}\right\Vert
_{L^{1}}^{1/2}\left\Vert q^{in}\right\Vert _{L^{\infty }}^{-1/2}$ in the 2D
vortex patch case. Contrary to this, if one is interested in bounded domains
or domains with periodic boundary conditions, typical length scale will be
dictated by the size of the domain, and then the rated of convergence will
be expressed as dimensionless quantities.

In the following all the constants $C$ are independent of $\alpha$, and all
the $\alpha $ dependencies are spelled explicitly.

\section{\label{sec:classical_solutions}Classical solutions}

In this section we study the convergence rate of the solutions of the Euler-$%
\alpha $ equations to the solution of the Euler equations for strong
solutions that belong to the Sobolev space $H^{m}\left( \mathbb{R}%
^{n}\right) $, $m>n/2+1$. For 2D incompressible flow with initial velocity $%
v^{in}\in H^{m}\left( \mathbb{R}^{2}\right) $, $m\geq 3$, the unique
solutions of the Euler equations exist globally in time (see, e.g., \cite%
{b_MB02}). Similarly to the Euler case, the Euler-$\alpha $ equations has a
unique global solution, since, as in the 2D Euler case, we have an \textit{a
priori} uniform control over the $L^{\infty }$ norm of the vorticity, which
implies the global existence, as in the proof of the Beale-Kato-Majda
criterion \cite{a_BKM84}. In $\mathbb{R}^{3}$ only local in time existence
of strong solutions of the Euler equations has been shown, see \cite%
{a_L25,a_L30,a_K72,a_K75,b_MB02}. The existence and well-posedness of the
Euler-$\alpha $ equations for a short time can be easily shown following the
classical theory of the Euler equations; see also \cite{a_HL07,a_LJ08} for
an analogue of the Beale-Kato-Majda criterion for the Euler-$\alpha $ model.
More precisely, one has the following result:

\begin{proposition}
\label{prop:strong_sol_existence}Let $v^{in}\in H^{m}\left( \mathbb{R}%
^{n}\right) $, $m>n/2+1$. There exists $T^{\ast }=T^{\ast }\left(
\left\Vert v^{in}\right\Vert _{H^{m}}\right) $, $T^{\ast }\geq
\frac{C}{\left\Vert v^{in}\right\Vert _{H^{m}}}$, such that for any
$T<T^{\ast }$ there exists a unique solution $v\in
C([0,T];H^{m}\left( \mathbb{R}^{n}\right) )\cap AC\left(
[0,T];H^{m-1}\left( \mathbb{R}^{n}\right) \right) $ of the Euler
equations \eqref{grp:EulerEq} with initial data $v^{in}$. In two
dimensions the solution exists globally in time. Similar results
hold for the Euler-$\alpha $ equations \eqref{grp:EulerAlphaEq} with
the maximal interval of existence of the three-dimensional
Euler-$\alpha $ equations being also dependent on $\alpha $.
\end{proposition}

In the next theorem we show that the solutions of the Euler-$\alpha $
equations for the $H^{m}$, $m>n/2+1$, initial data, exist at least as long
as the solution of the Euler system exists, and converge, as $\alpha
\rightarrow 0$, to the solution of Euler equations. Our proof of this result
follows the ideas in \cite{a_M07}.

\begin{theorem}
\label{thm:smoothCase}Let $v^{in},v^{in,\alpha }\in H^{m}\left( \mathbb{R}%
^{n}\right) $, $m>n/2+1$, and $\left\Vert v^{in}-v^{in,\alpha }\right\Vert
_{H^{m}}\rightarrow 0$, as $\alpha \rightarrow 0$. Let $T^{\ast }$ be the
time of existence of the solution of the Euler system \eqref{grp:EulerEq} $%
v\in C_{loc}([0,T^{\ast });H^{m}\left( \mathbb{R}^{n}\right) )\cap
AC_{loc}\left( [0,T^{\ast });H^{m-1}\left( \mathbb{R}^{n}\right) \right) $
with initial data $v^{in}$. Then, for all $0<T<T^{\ast }$ there exists $0<%
\bar{\alpha}=\bar{\alpha}\left( v^{in},v^{in,\alpha },T\right) $ such that
for all $\alpha \leq \bar{\alpha}$ there is a unique solution $v^{\alpha
}\in C([0,T];H^{m}\left( \mathbb{R}^{n}\right) )\cap AC\left(
[0,T];H^{m-1}\left( \mathbb{R}^{n}\right) \right) $ of the Euler-$\alpha $
equations \eqref{grp:EulerAlphaEq} with initial data $v^{in,\alpha }$.
Moreover,
\begin{equation*}
\left\Vert v^{\alpha }-v\right\Vert _{L^{\infty }\left( \left[ 0,T\right]
,H^{m}\right) }\rightarrow 0,
\end{equation*}%
as $\alpha \rightarrow 0$, and for all $0\leq t\leq T$
\begin{align}
& \left\Vert \left( v^{\alpha }-v\right) \left( t\right) \right\Vert
_{H^{m-2}}\leq \left( \left\Vert v^{in}-v^{in,\alpha }\right\Vert
_{H^{m-2}}+C\alpha ^{2}t\right) e^{Ct},  \label{eq:rate_s-2} \\
& \left\Vert \left( v^{\alpha }-v\right) \left( t\right) \right\Vert
_{H^{m-1}}\leq \left( \left\Vert v^{in}-v^{in,\alpha }\right\Vert
_{H^{m-1}}+C\alpha t\right) e^{Ct},  \label{eq:rate_s-1}
\end{align}%
where $C=C\left( \left\Vert v^{in}\right\Vert _{H^{m}},T\right) $ is
independent of $\alpha $.
\end{theorem}

We use the standard $H^{m}\left( \mathbb{R}^{n}\right) $ norm defined by%
\begin{equation*}
\left\Vert f\right\Vert _{H^{m}\left( \mathbb{R}^{n}\right) }^{2}=\int_{%
\mathbb{R}^{n}}\left( 1+\left\vert \xi \right\vert ^{2}\right)
^{m}\left\vert \hat{f}\left( \xi \right) \right\vert ^{2}d\xi ,
\end{equation*}%
where $\hat{f}$ denotes the Fourier transform of $f$.

To prove the theorem we need the following estimates (see, e.g., \cite%
{a_K72,b_MP94})

\begin{lemma}
\label{lem:non_lin_estimates}Let $n=2,3$. Let $u,v\in H^{m}\left( \mathbb{R}%
^{n}\right) ,\diver u=0$, then there exist a constant $C>0$,
depending on $m$, such that
\begin{equation}
\begin{array}{ll}
\left\vert \left( (u\cdot \nabla )v,v\right) _{H^{m}}\right\vert \leq
C\left\Vert \nabla u\right\Vert _{H^{m-1}}\left\Vert v\right\Vert
_{H^{m}}^{2}\leq C\left\Vert u\right\Vert _{H^{m}}\left\Vert v\right\Vert
_{H^{m}}^{2}, & m>\frac{n}{2}+1, \\
\left\vert \left( (u\cdot \nabla )v,v\right) _{H^{m}}\right\vert \leq
C\left\Vert \nabla u\right\Vert _{H^{2}}\left\Vert v\right\Vert
_{H^{m}}^{2}\leq C\left\Vert u\right\Vert _{H^{3}}\left\Vert v\right\Vert
_{H^{m}}^{2}, & m\leq \frac{n}{2}+1.%
\end{array}
\label{eq:non_lin_H_s_estimate1}
\end{equation}%
Let $u\in H^{m}\left( \mathbb{R}^{n}\right) $, $v\in H^{m+1}\left( \mathbb{R}%
^{n}\right) $, and let $\Psi \left( u,v\right) $ be either one of the
following bilinear forms: $\Psi \left( u,v\right) =$ $u\times \left( \nabla
\times v\right) $, $\Psi \left( u,v\right) =$ $(u\cdot \nabla )v$ or $\Psi
\left( u,v\right) =v_{_{j}}\nabla {u}_{j}$, then there exists a constant $C>0
$, depending on $m$, such that
\begin{equation}
\begin{array}{ll}
\left\Vert \Psi \left( u,v\right) \right\Vert _{H^{m}}\leq C\left\Vert
u\right\Vert _{H^{m}}\left\Vert v\right\Vert _{H^{m+1}}, & m>\frac{n}{2}, \\
\left\Vert \Psi \left( u,v\right) \right\Vert _{H^{m}}\leq C\left\Vert
u\right\Vert _{H^{m}}\left\Vert v\right\Vert _{H^{3}}, & m\leq \frac{n}{2},%
\end{array}
\label{eq:non_lin_H_s_estimate2}
\end{equation}%
and for $m>\frac{n}{2}+1$
\begin{equation}
\left\Vert \Psi \left( u,v\right) \right\Vert _{H^{m}}\leq C\left(
\left\Vert u\right\Vert _{H^{m}}\left\Vert v\right\Vert _{H^{m}}+\left\Vert
u\right\Vert _{L^{\infty }}\left\Vert v\right\Vert _{H^{m+1}}\right) .
\label{eq:non_lin_H_s_estimate3}
\end{equation}
\end{lemma}

We also use the following lemma

\begin{lemma}
Let $g\in L^{2}\left( \mathbb{R}^{n}\right) $ and $f=\left( 1-\alpha
^{2}\Delta \right) ^{-1}g$. Then
\begin{equation}
\alpha \left\Vert \left( -\Delta \right) ^{1/2}f\right\Vert _{L^{2}}\leq
\frac{1}{2}\left\Vert g\right\Vert _{L^{2}}.  \label{eq:alpha_L_2_grad_u}
\end{equation}
\end{lemma}

\begin{proof}
By taking a Fourier transform we have
\begin{equation*}
\hat{f}\left( \xi \right) =\frac{\hat{g}\left( \xi \right) }{\left( 1+\alpha
^{2}\left\vert \xi \right\vert ^{2}\right) },
\end{equation*}%
hence%
\begin{align*}
\alpha ^{2}\left\Vert \left( -\Delta \right) ^{1/2}f\right\Vert
_{L^{2}\left( \mathbb{R}^{n}\right) }^{2}& =\int_{\mathbb{R}^{n}}\left\vert
\hat{g}\left( \xi \right) \right\vert ^{2}\frac{\alpha ^{2}\left\vert \xi
\right\vert ^{2}}{\left( 1+\alpha ^{2}\left\vert \xi \right\vert ^{2}\right)
^{2}}d\xi  \\
& \leq \left\Vert g\right\Vert _{L^{2}\left( \mathbb{R}^{n}\right)
}^{2}\sup_{y\geq 0}\frac{y}{\left( 1+y\right) ^{2}}\leq \frac{1}{4}%
\left\Vert g\right\Vert _{L^{2}\left( \mathbb{R}^{n}\right) }^{2}.
\end{align*}
\end{proof}

Next, we prove Theorem \ref{thm:smoothCase}.

\begin{proof}
First, let us assume that $v,v^{\alpha }\in C([0,T];H^{m}\left( \mathbb{R}%
^{n}\right) )\cap AC\left( [0,T];H^{m-1}\left( \mathbb{R}^{n}\right) \right)
$ are the solutions of the Euler and the Euler-$\alpha $ systems %
\eqref{grp:EulerEq} and \eqref{grp:EulerAlphaEq} with initial data $v^{in}$
and $v^{in,\alpha }$, respectively, on some mutual time interval $[0,T]$,
with
\begin{equation}
\left\Vert v\right\Vert _{L^{\infty }\left( \left[ 0,T\right] ,H^{m}\right) }%
\text{,}\left\Vert v^{\alpha }\right\Vert _{L^{\infty }\left( \left[ 0,T%
\right] ,H^{m}\right) }\leq C\left( T,\left\Vert v^{in}\right\Vert
_{H^{m}}\right) ,  \label{eq:bnd_on_H_m_of_v_and_v_alpha}
\end{equation}%
and let us show the convergence rates in $H^{m-2}$ and $H^{m-1}$. In the
second part we show that the solutions of the Euler-$\alpha $ equations
exist at least up to the time of existence of the solution of Euler
equations and satisfy \eqref{eq:bnd_on_H_m_of_v_and_v_alpha}.

The difference $w^{\alpha }=v-v^{\alpha }$ satisfies the following equation%
\begin{align}
& \frac{\partial w^{\alpha }}{\partial t}+(w^{\alpha }\cdot \nabla
)v+v_{j}\nabla w_{j}^{\alpha }+({u}^{\alpha }\cdot \nabla )w^{\alpha
}+w_{_{j}}^{\alpha }\nabla {u}_{j}^{\alpha }  \label{eq:diff_eq_smooth_case}
\\
& \qquad -(\alpha ^{2}\Delta {u}^{\alpha }\cdot \nabla )v-v_{j}\nabla \left(
\alpha ^{2}\Delta {u}_{j}^{\alpha }\right) -v_{j}\nabla v_{j}+\nabla \left(
p-p^{\alpha }\right) =0,  \notag \\
& w^{\alpha }(x,0)=v^{in}(x,0)-v^{in,\alpha }(x,0).  \notag
\end{align}%
Let $k$ be either $m-1$ or $m-2$. Taking the $H^{k}$ inner product of %
\eqref{eq:diff_eq_smooth_case} with $w^{\alpha }\left( t\right) $, we obtain
\begin{equation*}
\frac{1}{2}\frac{d}{dt}\left\Vert w^{\alpha }\right\Vert _{H^{k}}^{2}\leq
I_{1}+I_{2}+I_{3}+I_{4},
\end{equation*}
where%
\begin{align*}
I_{1}& =\left\vert \left( w^{\alpha }\times \left( \nabla \times v\right)
,w^{\alpha }\right) _{H^{k}}\right\vert , \\
I_{2}& =\left\vert \left( ({u}^{\alpha }\cdot \nabla )w^{\alpha },w^{\alpha
}\right) _{H^{k}}\right\vert , \\
I_{3}& =\left\vert \left( w_{_{j}}^{\alpha }\nabla {u}_{j}^{\alpha
},w^{\alpha }\right) _{H^{k}}\right\vert , \\
I_{4}& =\alpha ^{2}\left\vert \left( \Delta {u}^{\alpha }\times \left(
\nabla \times v\right) ,w^{\alpha }\right) _{H^{k}}\right\vert ,
\end{align*}%
due to the identity
\begin{equation}
\left( b\cdot \nabla \right) a+a_{j}\nabla b_{j}=-b\times \left( \nabla
\times a\right) +\nabla \left( a\cdot b\right) .  \label{eq:vector_identity}
\end{equation}%
Now, by \eqref{eq:non_lin_H_s_estimate2}%
\begin{equation*}
I_{1}\leq C\left\Vert v\right\Vert _{H^{m}}\left\Vert w^{\alpha }\right\Vert
_{H^{k}}^{2},
\end{equation*}%
by \eqref{eq:non_lin_H_s_estimate1}%
\begin{equation*}
I_{2}\leq C\left\Vert {u}^{\alpha }\right\Vert _{H^{m}}\left\Vert w^{\alpha
}\right\Vert _{H^{k}}^{2}\leq C\left\Vert {v}^{\alpha }\right\Vert
_{H^{m}}\left\Vert w^{\alpha }\right\Vert _{H^{k}}^{2},
\end{equation*}%
by \eqref{eq:non_lin_H_s_estimate2}%
\begin{equation*}
I_{3}\leq C\left\Vert {u}^{\alpha }\right\Vert _{H^{m}}\left\Vert w^{\alpha
}\right\Vert _{H^{k}}^{2}\leq C\left\Vert {v}^{\alpha }\right\Vert
_{H^{m}}\left\Vert w^{\alpha }\right\Vert _{H^{k}}^{2},
\end{equation*}%
and by \eqref{eq:non_lin_H_s_estimate2}%
\begin{equation*}
I_{4}\leq \alpha ^{2}\left\Vert \Delta {u}^{\alpha }\right\Vert
_{H^{k}}\left\Vert {v}\right\Vert _{H^{m}}\left\Vert w^{\alpha }\right\Vert
_{H^{k}}.
\end{equation*}%
Summing up we get the energy estimate%
\begin{equation*}
\frac{1}{2}\frac{d}{dt}\left\Vert w^{\alpha }\right\Vert _{H^{k}}^{2}\leq
C\left\Vert w^{\alpha }\right\Vert _{H^{k}}^{2}\left( \left\Vert
v\right\Vert _{H^{m}}+\left\Vert {v}^{\alpha }\right\Vert _{H^{m}}\right)
+\alpha ^{2}\left\Vert \Delta {u}^{\alpha }\right\Vert _{H^{k}}\left\Vert {v}%
\right\Vert _{H^{m}}\left\Vert w^{\alpha }\right\Vert _{H^{k}}.
\end{equation*}
Now, for $k=m-2$ we have that
\begin{equation*}
\left\Vert \Delta {u}^{\alpha }\right\Vert _{H^{m-2}}\leq C\left\Vert {u}%
^{\alpha }\right\Vert _{H^{m}}\leq C\left\Vert {v}^{\alpha }\right\Vert
_{H^{m}},
\end{equation*}%
and hence%
\begin{equation*}
\frac{1}{2}\frac{d}{dt}\left\Vert w^{\alpha }\right\Vert _{H^{m-2}}^{2}\leq
C\left\Vert w^{\alpha }\right\Vert _{H^{m-2}}^{2}\left( \left\Vert
v\right\Vert _{H^{m}}+\left\Vert {v}^{\alpha }\right\Vert _{H^{m}}\right)
+C\alpha ^{2}\left\Vert {v}^{\alpha }\right\Vert _{H^{m}}\left\Vert {v}%
\right\Vert _{H^{m}}\left\Vert w^{\alpha }\right\Vert _{H^{m-2}}.
\end{equation*}%
while for $k=m-1$ we have by \eqref{eq:alpha_L_2_grad_u}
\begin{equation*}
\alpha \left\Vert \Delta {u}^{\alpha }\right\Vert _{H^{m-1}}\leq C\left\Vert
\left( -\Delta \right) ^{1/2}{v}^{\alpha }\right\Vert _{H^{m-1}}\leq
C\left\Vert v^{\alpha }\right\Vert _{H^{m}},
\end{equation*}%
and hence%
\begin{equation*}
\frac{1}{2}\frac{d}{dt}\left\Vert w^{\alpha }\right\Vert _{H^{m-1}}^{2}\leq
C\left\Vert w^{\alpha }\right\Vert _{H^{m-1}}^{2}\left( \left\Vert
v\right\Vert _{H^{m}}+\left\Vert {v}^{\alpha }\right\Vert _{H^{m}}\right)
+C\alpha \left\Vert v^{\alpha }\right\Vert _{H^{m}}\left\Vert {v}\right\Vert
_{H^{m}}\left\Vert w^{\alpha }\right\Vert _{H^{m-1}}.
\end{equation*}%
Therefore, by Gr\"{o}nwall lemma and \eqref{eq:bnd_on_H_m_of_v_and_v_alpha}
we obtain that%
\begin{equation*}
\left\Vert w^{\alpha }\left( \cdot ,t\right) \right\Vert _{H^{m-2}}\leq
\left( \left\Vert w^{\alpha }\left( \cdot ,0\right) \right\Vert
_{H^{m-2}}+C\alpha ^{2}t\right) e^{Ct},
\end{equation*}%
and%
\begin{equation*}
\left\Vert w^{\alpha }\left( \cdot ,t\right) \right\Vert _{H^{m-1}}\leq
\left( \left\Vert w^{\alpha }\left( \cdot ,0\right) \right\Vert
_{H^{m-1}}+C\alpha t\right) e^{Ct},
\end{equation*}%
where $C=C\left( \left\Vert v^{in}\right\Vert _{H^{m}},T\right) $.

Next, we show that the solution of the Euler-$\alpha $ equations %
\eqref{grp:EulerAlphaEq} exists as long as we can solve the Euler system %
\eqref{grp:EulerEq}, and that $v^{\alpha }$ converges to $v$ in $L^{\infty
}\left( \left[ 0,T\right] ,H^{m}\right) $, as $\alpha \rightarrow 0$. The
proof follows the ideas in \cite{a_M07}. We regularize the initial data by
taking $v^{in,\delta }=\mathcal{F}^{-1}\left( \chi _{\left\vert \xi
\right\vert \leq 1/\delta }\left( \xi \right) \mathcal{F}\left(
v^{in}\right) \right) $ and $v^{in,\delta ,\alpha }=\mathcal{F}^{-1}\left(
\chi _{\left\vert \xi \right\vert \leq 1/\delta }\left( \xi \right) \mathcal{%
F}\left( v^{in,\alpha }\right) \right) $, for some $\delta \in \left(
0,\delta _{0}\right] $. Let $\alpha ^{\ast }$ be such that $v^{in,\alpha
^{\ast }}\neq 0$ (if $v^{in,\alpha }=0$ for all $\alpha $, then the proof is
trivial). Since $\left\Vert v^{in,\alpha }-v^{in}\right\Vert
_{H^{m}}\rightarrow 0$, as $\alpha \rightarrow 0$, then there exists $\alpha
_{0}$ such that $\left\Vert v^{in,\alpha }\right\Vert _{H^{m}}\leq
\left\Vert v^{in}\right\Vert _{H^{m}}+\left\Vert v^{in,\alpha ^{\ast
}}\right\Vert _{H^{m}}$ for all $\alpha \leq \alpha _{0}$. The regularized
initial velocities of the Euler and of the Euler-$\alpha $ equations
satisfy, for all $\alpha \leq \alpha _{0}$,
\begin{align*}
& \left\Vert v^{in,\alpha ,\delta }\right\Vert _{H^{m}},\left\Vert
v^{in,\delta }\right\Vert _{H^{m}}\leq K, \\
& \left\Vert v^{in,\alpha ,\delta }\right\Vert _{H^{m+1}},\left\Vert
v^{in,\delta }\right\Vert _{H^{m+1}}\leq \frac{K}{\delta },
\end{align*}%
and for $s\in \left[ 0,m\right] +$
\begin{equation*}
\left\Vert v^{in,\alpha }-v^{in,\alpha ,\delta }\right\Vert
_{H^{s}},\left\Vert v^{in}-v^{in,\delta }\right\Vert _{H^{s}}\leq K\delta
^{m-s},
\end{equation*}%
where $K=K\left( \left\Vert v^{in}\right\Vert _{H^{m}},\left\Vert
v^{in,\alpha ^{\ast }}\right\Vert _{H^{m}}\right) $. In the following we fix
$s$ such that $\frac{n}{2}<s<m-1$.

Let $v$, $v^{\delta }$, $v^{\alpha }$, $v^{\alpha ,\delta }$ be the
corresponding solutions of Euler and Euler-$\alpha $ equations with initial
data $v^{in}$ and $v^{in,\delta }$, respectively. Notice that, as explained
below, the solutions $v$, $v^{\delta }$, $v^{\alpha }$ and $v^{\alpha
,\delta }$ exist on some time interval $\left[ 0,T_{0}\right] $, $\frac{C}{K}%
\leq T_{0}<T^{\ast }$, which is independent of $\alpha $ and $\delta $, and
also, for all $t\in \left[ 0,T_{0}\right] $%
\begin{equation}
\left\Vert \varphi \left( t\right) \right\Vert _{H^{m}},\left\Vert \varphi
^{\delta }\left( t\right) \right\Vert _{H^{m}}\leq C\left( T_{0},K\right)
,\;\left\Vert \varphi ^{\delta }\left( t\right) \right\Vert _{H^{m+1}}\leq
\frac{C\left( T_{0},K\right) }{\delta },  \label{eq:sol_bnd}
\end{equation}%
where $\varphi $ denotes $v$ or $v^{\alpha }$, and $\varphi ^{\delta }$
denotes either $v^{\delta }$ or $v^{\alpha ,\delta }$. Indeed, in $\mathbb{R}%
^{3}$, writing Euler and Euler-$\alpha $ equations in the vorticity
formulation we have that $q=\curl v$ and $q^{\alpha }=\curl v^{\alpha }$
satisfy
\begin{equation*}
\frac{\partial q}{\partial t}+\left( v\cdot \nabla \right) q=\left( q\cdot
\nabla \right) v
\end{equation*}%
and
\begin{equation*}
\frac{\partial q^{\alpha }}{\partial t}+\left( u^{\alpha }\cdot \nabla
\right) q^{\alpha }=\left( q^{\alpha }\cdot \nabla \right) u^{\alpha },
\end{equation*}%
respectively. Making the $H^{k}$ estimates (one can use %
\eqref{eq:non_lin_H_s_estimate1} and the fact that for $k>3/2$, $H^{k}\left(
\mathbb{R}^{3}\right) $ is a Banach algebra), due to $m-1>\frac{3}{2}$, we
obtain%
\begin{align*}
& \frac{d}{dt}\left\Vert \psi \right\Vert _{H^{m-1}}\leq C\left\Vert \psi
\right\Vert _{H^{m-1}}^{2}, \\
& \frac{d}{dt}\left\Vert \psi ^{\delta }\right\Vert _{H^{m-1}}\leq
C\left\Vert \psi ^{\delta }\right\Vert _{H^{m-1}}^{2}, \\
& \frac{d}{dt}\left\Vert \psi ^{\delta }\right\Vert _{H^{m}}\leq C\left\Vert
\psi ^{\delta }\right\Vert _{H^{m-1}}\left\Vert \psi ^{\delta }\right\Vert
_{H^{m}},
\end{align*}%
where $\psi $ denotes $q$ or $q^{\alpha }$, and $\psi ^{\delta }$ denotes $%
q^{\delta }$ or $q^{\alpha ,\delta }$, and hence the solutions exist on a
certain interval $\left[ 0,T_{0}\right] $ depending only on $K$ and $\delta
_{0}$, independent of $\alpha $ and $\delta $. Furthermore, we have\footnote{%
It follows that
\begin{equation}
\left\Vert \psi \left( t\right) \right\Vert _{H^{m-1}}\leq \Psi \left(
t\right) ,  \label{eq:ft1}
\end{equation}%
where $\Psi $ is the solution of the scalar value initial-value problem $%
\frac{d}{dt}\Psi \left( t\right) =C\Psi ^{2}\left( t\right) $, where $\Psi
\left( 0\right) $ is the $H^{m-1}$ norm of either one of $q^{in}$, $%
q^{in,\alpha }$, $q^{in,\delta }$, $q^{in,\alpha ,\delta }$, $\Psi \left(
0\right) \leq K$.
\par
Now $\Psi \left( t\right) =\frac{\Psi \left( 0\right) }{1-Ct\Psi \left(
0\right) }$ exists on a certain interval of time $\left[ 0,T_{0}\right] $, $%
T_{0}>\frac{1}{C\Psi \left( 0\right) }\geq \frac{1}{CK}$, obviously
independent of $\alpha $ and $\delta $, and we have $\sup_{0\leq t\leq
T_{0}}\Psi \left( t\right) \leq \frac{K}{1-CT_{0}K}$. Now, if a solution $%
\psi \left( t\right) $ existing on $\left[ 0,T^{\alpha ,\delta }\right] $,
such that $T^{\alpha ,\delta }<T_{0}$, then the system can be solved with
initial value $\psi ^{\delta }\left( T^{\alpha ,\delta }\right) \in H^{m-1}$%
, to continue the solution to $\left[ 0,T^{\alpha ,\delta }+T_{1}^{\alpha
,\delta }\right] $ in which \eqref{eq:ft1} is true. Iterating this argument
we can continue the solution to cover the whole interval $\left[ 0,T\right] $
with the estimate \eqref{eq:ft1} throughout. (If the solution cannot be
continued at some time $\tilde{T}<T_{0}$, then $\lim \sup_{t\rightarrow
\tilde{T}^{-}}\left\Vert \psi ^{\delta }\left( t\right) \right\Vert
_{H^{m-1}}=\infty $, a contradiction to \eqref{eq:ft1}.)}%
\begin{equation*}
\left\Vert \psi \left( t\right) \right\Vert _{H^{m-1}}\leq C\left(
T_{0},K\right) ,\left\Vert \psi ^{\delta }\left( t\right) \right\Vert
_{H^{m}}\leq \frac{C\left( T_{0},K\right) }{\delta }.
\end{equation*}%
Hence, from the Biot-Savart law, $\varphi =\frac{1}{4\pi }\int_{\mathbb{R}%
^{3}}\frac{\left( x-y\right) }{\left\vert x-y\right\vert ^{3}}\times \curl%
\varphi \left( y\right) dy$, we obtain \eqref{eq:sol_bnd}. In $\mathbb{R}^{2}
$ the solutions of Euler and Euler-$\alpha $ equations exist for all times,
see Proposition \ref{prop:strong_sol_existence}.

Now, one can show, following \cite{a_M07}, that $v^{\delta }$ is a Cauchy
sequence in the Banach space $C\left( \left[ 0,T_{0}\right] ,H^{m}\right) $,
due to
\begin{equation}
\left\Vert v^{\delta }-v^{\delta ^{\prime }}\right\Vert _{_{C\left( \left[
0,T_{0}\right] ,H^{m}\right) }}\leq \left( \left\Vert v^{in,\delta ^{\prime
}}-v^{in,\delta }\right\Vert _{H^{m}}+\left( \max \left\{ \delta ,\delta
^{\prime }\right\} \right) ^{m-s-1}C\left( T_{0},\left\Vert
v^{in}\right\Vert _{H^{m}}\right) \right) e^{C\left( T_{0},\left\Vert
v^{in}\right\Vert _{H^{m}}\right) }.
\label{eq:convergence_reg_Euler_to_Euler_CauchyEst}
\end{equation}%
The limit $v$ of $v^{\delta }$,as $\delta \rightarrow 0$, is the solution of
Euler equations, and from \eqref{eq:convergence_reg_Euler_to_Euler_CauchyEst}%
, we have
\begin{equation}
\left\Vert v-v^{\delta }\right\Vert _{L^{\infty }\left( \left[ 0,T_{0}\right]
,H^{m}\right) }\leq \left( \left\Vert v^{in}-v^{in,\delta }\right\Vert
_{H^{m}}+C\left( T_{0},\left\Vert v^{in}\right\Vert _{H^{m}}\right) \delta
^{m-s-1}\right) Ce^{C\left( T_{0},\left\Vert v^{in}\right\Vert
_{H^{m}}\right) }.  \label{eq:convergence_reg_Euler_to_Euler}
\end{equation}

Next we show that $v^{\alpha ,\delta }$ is a Cauchy sequence in $C\left( %
\left[ 0,T_{0}\right] ,H^{m}\right) $, and converges to the solution of
Euler-$\alpha $ equations $v^{\alpha }$, as $\delta \rightarrow 0$, and we
have
\begin{equation}
\left\Vert v^{\alpha }-v^{\alpha ,\delta }\right\Vert _{L^{\infty }\left(
\left[ 0,T_{0}\right] ,H^{m}\right) }\leq \left( \left\Vert v^{in,\alpha
}-v^{in,\alpha ,\delta }\right\Vert _{H^{m}}+C\left( T_{0},K\right) \delta
^{m-s-1}\right) Ce^{C\left( T_{0},K\right) }.
\label{eq:convergence_reg_Euler_alpha_to_Euler_alpha}
\end{equation}%
Specifically, in $\mathbb{R}^{3}$ ($\mathbb{R}^{2}$ case can be done
similarly), we first assume that $\delta >\delta ^{\prime }$, and denote $%
\varpi ^{\alpha ,\delta ,\delta ^{\prime }}=q^{\alpha ,\delta ^{\prime
}}-q^{\alpha ,\delta }$, $w^{\alpha ,\delta ,\delta ^{\prime }}=v^{\alpha
,\delta ^{\prime }}-v^{\alpha ,\delta }$. We have that%
\begin{equation*}
\frac{\partial \varpi ^{\alpha ,\delta ,\delta ^{\prime }}}{\partial t}%
+\left( u^{\alpha ,\delta ^{\prime }}\cdot \nabla \right) \varpi ^{\alpha
,\delta ,\delta ^{\prime }}+\left( \left( \left( 1-\alpha ^{2}\Delta \right)
^{-1}w^{\alpha ,\delta ,\delta ^{\prime }}\right) \cdot \nabla \right)
q^{\alpha ,\delta }=\left( \varpi ^{\alpha ,\delta ,\delta ^{\prime }}\cdot
\nabla \right) u^{\alpha ,\delta ^{\prime }}+\left( q^{\alpha ,\delta }\cdot
\nabla \right) \left( \left( 1-\alpha ^{2}\Delta \right) ^{-1}w^{\alpha
,\delta ,\delta ^{\prime }}\right) .
\end{equation*}%
Using \eqref{eq:non_lin_H_s_estimate1} and \eqref{eq:non_lin_H_s_estimate3}
we obtain
\begin{equation*}
\frac{d}{dt}\left\Vert \varpi ^{\alpha ,\delta ,\delta ^{\prime
}}\right\Vert _{H^{m-1}}^{2}\leq C\left( \left\Vert q^{\alpha ,\delta
^{\prime }}\right\Vert _{H^{m-1}}+\left\Vert q^{\alpha ,\delta }\right\Vert
_{H^{m-1}}\right) \left\Vert \varpi ^{\alpha ,\delta ,\delta ^{\prime
}}\right\Vert _{H^{m-1}}^{2}+C\left\Vert w^{\alpha ,\delta ,\delta ^{\prime
}}\right\Vert _{L^{\infty }}\left\Vert q^{\alpha ,\delta }\right\Vert
_{H^{m}}\left\Vert \varpi ^{\alpha ,\delta ,\delta ^{\prime }}\right\Vert
_{H^{m-1}}.
\end{equation*}%
Now, the difference $w^{\alpha ,\delta ,\delta ^{\prime }}=v^{\alpha ,\delta
^{\prime }}-v^{\alpha ,\delta }$ satisfies the equation (see %
\eqref{eq:vector_identity})%
\begin{equation*}
\frac{\partial w^{\alpha ,\delta ,\delta ^{\prime }}}{\partial t}+\left(
\left( 1-\alpha ^{2}\Delta \right) ^{-1}w^{\alpha ,\delta }\right) \times
\left( \nabla \times {v}^{\alpha ,\delta ^{\prime }}\right) +u^{\alpha
,\delta }\times \left( \nabla \times w^{\alpha ,\delta ,\delta ^{\prime
}}\right) +\nabla \left( \tilde{p}^{\alpha ,\delta ^{\prime }}-\tilde{p}%
^{\alpha ,\delta }\right) =0,
\end{equation*}%
where $\tilde{p}^{\alpha }$, $\tilde{p}^{\alpha ,\delta }$ are the modified
pressure. By \eqref{eq:non_lin_H_s_estimate1} and %
\eqref{eq:non_lin_H_s_estimate2} we obtain that%
\begin{equation*}
\frac{1}{2}\frac{d}{dt}\left\Vert w^{\alpha ,\delta ,\delta ^{\prime
}}\right\Vert _{H^{s}}^{2}\leq C\left( \left\Vert {v}^{\alpha ,\delta
}\right\Vert _{H^{m}}+\left\Vert {v}^{\alpha ,\delta ^{\prime }}\right\Vert
_{H^{m}}\right) \left\Vert w^{\alpha ,\delta ,\delta ^{\prime }}\right\Vert
_{H^{s}}^{2}.
\end{equation*}%
Recall, that $\frac{n}{2}<s$, hence, by Sobolev embedding theorem and Gr\"{o}%
nwall lemma we have%
\begin{eqnarray*}
\left\Vert w^{\alpha ,\delta ,\delta ^{\prime }}\right\Vert _{L^{\infty
}\left( \left[ 0,T_{0}\right] ;L^{\infty }\right) } &\leq &\left\Vert
w^{\alpha ,\delta ,\delta ^{\prime }}\right\Vert _{L^{\infty }\left( \left[
0,T\right] ;H^{s}\right) } \\
&\leq &\left\Vert w^{\alpha ,\delta ,\delta ^{\prime }}\left( 0\right)
\right\Vert _{H^{s}}e^{C\int_{0}^{t}\left\Vert {v}^{\alpha ,\delta }\left(
\tau \right) \right\Vert _{H^{m}}+\left\Vert {v}^{\alpha ,\delta ^{\prime
}}\left( \tau \right) \right\Vert _{H^{m}}d\tau } \\
&\leq &C\left( T_{0},K\right) \delta ^{m-s}e^{C\left( T_{0},K\right) }.
\end{eqnarray*}%
Using also that $\left\Vert {q}^{\alpha ,\delta }\right\Vert _{H^{m}}\leq
\frac{C\left( T_{0},K\right) }{\delta }$ we obtain%
\begin{equation*}
\frac{d}{dt}\left\Vert \varpi ^{\alpha ,\delta ,\delta ^{\prime
}}\right\Vert _{H^{m-1}}\leq C\left( \left\Vert q^{\alpha ,\delta ^{\prime
}}\right\Vert _{H^{m-1}}+\left\Vert q^{\alpha ,\delta }\right\Vert
_{H^{m-1}}\right) \left\Vert \varpi ^{\alpha ,\delta ,\delta ^{\prime
}}\right\Vert _{H^{m-1}}+C\left( T_{0},K\right) \delta ^{m-s-1}e^{C\left(
T_{0},K\right) }
\end{equation*}%
hence
\begin{equation*}
\left\Vert \varpi ^{\alpha ,\delta }\right\Vert _{L^{\infty }\left( \left[
0,T_{0}\right] ;H^{m-1}\right) }\leq \left( \left\Vert q^{in,\alpha ,\delta
^{\prime }}-q^{in,\alpha ,\delta }\right\Vert _{H^{m-1}}+\delta
^{m-s-1}C\left( T_{0},K\right) \right) e^{C\left( T_{0},K\right) }.
\end{equation*}%
Making the same estimates for the case $\delta ^{\prime }\geq \delta
$, we obtain that
\begin{equation*}
\left\Vert \varpi ^{\alpha ,\delta }\right\Vert _{L^{\infty }\left( \left[
0,T_{0}\right] ;H^{m-1}\right) }\leq \left( \left\Vert q^{in,\alpha ,\delta
^{\prime }}-q^{in,\alpha ,\delta }\right\Vert _{H^{m-1}}+\left( \max \left\{
\delta ,\delta ^{\prime }\right\} \right) ^{m-s-1}C\left( T_{0},K\right)
\right) e^{C\left( T_{0},K\right) },
\end{equation*}%
which implies that $\varpi ^{\alpha ,\delta }$ is a Cauchy sequence in the
Banach space $C\left( \left[ 0,T_{0}\right] ,H^{m}\right) $, and its limit $%
v^{\alpha }$, which is the solution of Euler-$\alpha $ equations satisfies %
\eqref{eq:convergence_reg_Euler_alpha_to_Euler_alpha}.

It remains to show that $\left\Vert v^{\delta }\left( t\right) -v^{\alpha
,\delta }\left( t\right) \right\Vert _{H^{m}}$ converge to zero, as $\alpha $
and $\delta $ converge to zero uniformly in $\left[ 0,T_{0}\right] $, and
then the convergence of $v^{\alpha }$ to $v$ in $L^{\infty }\left( \left[
0,T_{0}\right] ,H^{m}\right) $ follows.

From \eqref{eq:diff_eq_smooth_case} and Lemma \ref{lem:non_lin_estimates} we
have that the difference $w^{\delta }=v^{\delta }-{v}^{\alpha ,\delta }$
satisfy%
\begin{align*}
\frac{d}{dt}\left\Vert w^{\delta }\right\Vert _{H^{m}}& \leq C\left(
\left\Vert v^{\delta }\right\Vert _{H^{m}}+\left\Vert {v}^{\alpha ,\delta
}\right\Vert _{H^{m}}\right) \left\Vert w^{\delta }\right\Vert
_{H^{m}}+C\left\Vert w^{\delta }\right\Vert _{L^{\infty }}\left( \left\Vert
v^{\delta }\right\Vert _{H^{m+1}}+\left\Vert {v}^{\alpha ,\delta
}\right\Vert _{H^{m+1}}\right) \\
& +C\alpha ^{2}\left\Vert \Delta {u}^{\alpha ,\delta }\right\Vert
_{H^{m}}\left\Vert v^{\delta }\right\Vert _{H^{m}}+C\alpha ^{2}\left\Vert
\Delta {u}^{\alpha ,\delta }\right\Vert _{L^{\infty }}\left\Vert v^{\delta
}\right\Vert _{H^{m+1}}.
\end{align*}%
By Sobolev embedding theorem, \eqref{eq:rate_s-1} and Gr\"{o}nwall lemma we
have that%
\begin{align*}
\left\Vert w^{\delta }\right\Vert _{L^{\infty }\left( \left[ 0,T_{0}\right]
,L^{\infty }\left( \mathbb{R}^{n}\right) \right) }& \leq \left\Vert
w^{\delta }\right\Vert _{L^{\infty }\left( \left[ 0,T_{0}\right]
,H^{m-1}\right) }\leq \left( \left\Vert v^{in,\delta }-v^{in,\alpha ,\delta
}\right\Vert _{H^{m-1}}+C\left( T_{0},K\right) \alpha \right) e^{C\left(
T_{0},K\right) } \\
& \leq \left( \left\Vert v^{in}-v^{in,\alpha }\right\Vert _{H^{m-1}}+C\left(
T_{0},K\right) \alpha \right) e^{C\left( T_{0},K\right) },
\end{align*}%
by \eqref{eq:alpha_L_2_grad_u}
\begin{align*}
& \alpha \left\Vert \Delta {u}^{\alpha }\right\Vert _{H^{m}}\leq C\left\Vert
v^{\alpha }\right\Vert _{H^{m+1}}, \\
& \alpha \left\Vert \Delta {u}^{\alpha ,\delta }\right\Vert _{L^{\infty
}}\leq \alpha \left\Vert \Delta {u}^{\alpha ,\delta }\right\Vert
_{H^{m-1}}\leq C\left\Vert v^{\alpha }\right\Vert _{H^{m}}.
\end{align*}%
Hence%
\begin{equation*}
\frac{d}{dt}\left\Vert w^{\delta }\right\Vert _{H^{m}}\leq C\left\Vert
w^{\delta }\right\Vert _{H^{m}}+C\left( \frac{\left\Vert v^{in}-v^{in,\alpha
}\right\Vert _{H^{m-1}}}{\delta }+\frac{\alpha }{\delta }\right) e^{C}+C%
\frac{\alpha }{\delta }
\end{equation*}%
where $C=C\left( K,T_{0}\right) $. By using Gr\"{o}nwall lemma, first
letting $\alpha \rightarrow 0$, and then letting $\delta \rightarrow 0$,
while choosing $\delta $ such that $\frac{\left\Vert v^{in}-v^{in,\alpha
}\right\Vert _{H^{m-1}}}{\delta }\rightarrow 0$, we obtain that $\left\Vert
v^{\delta }-v^{\alpha ,\delta }\right\Vert _{L^{\infty }\left( \left[ 0,T_{0}%
\right] ,H^{m}\right) }\rightarrow 0$. Using also %
\eqref{eq:convergence_reg_Euler_to_Euler} and %
\eqref{eq:convergence_reg_Euler_alpha_to_Euler_alpha} the convergence of $%
v^{\alpha }$ to $v$ in $L^{\infty }\left( \left[ 0,T_{0}\right]
,H^{m}\right) $ follows.

Now, let the Euler solution exist on $[0,T^{\ast })$, then we can continue
the solution of the Euler-$\alpha $ equations up to any $T<T^{\ast }$ in a
finite number of iterations using the above argument. Indeed, we continue
the solution from a time interval $\left[ 0,T_{0}+\ldots +T_{k-1}\right] $
to a time interval $\left[ 0,T_{0}+\ldots +T_{k}\right] $, by solving the
Euler-$\alpha $ equation with the initial data $v^{\alpha }\left(
T_{0}+\ldots +T_{k-1}\right) $, $\alpha \leq \min \left\{ \alpha _{0},\ldots
,\alpha _{k-1}\right\} $, which converges to $v\left( T_{0}+\ldots
+T_{k-1}\right) $ for a time $T_{k}\geq C\left( \left\Vert v\left(
T_{k-1}\right) \right\Vert _{H^{m}}+\left\Vert v^{in,\alpha ^{\ast
}}\right\Vert _{H^{m}}\right) ^{-1}\geq C\left( \left\Vert v\right\Vert
_{L^{\infty }\left( \left[ 0,T\right] ,H^{m}\right) }+\left\Vert
v^{in,\alpha ^{\ast }}\right\Vert _{H^{m}}\right) ^{-1}$.
\end{proof}

\section{\label{sec:vortex_patch} The vortex patch case}

In this section we study the convergence, as $\alpha \rightarrow 0$, of the $%
L^{2}$ norm of the difference between the solutions of the 2D Euler-$\alpha $
equations \eqref{grp:EulerAlphaEq} and the solution of the 2D Euler
equations \eqref{grp:EulerEq} for the vortex patch initial data with
vorticity being a characteristic function of a simply connected bounded
domain with $C^{1,\gamma }$, $\gamma \in \left( 0,1\right) $, boundary. We
show that the convergence rate is of order $\left( \alpha ^{2}\right) ^{3/4}$%
, which corresponds to the optimal convergence rate of the difference
between the solutions of the 2D NSE and the 2D Euler equation which is of
order $\left( \nu t\right) ^{3/4}$.

Yudovich \cite{a_Y63} obtained the existence and uniqueness of weak
solutions of Euler equations for initial vorticity in $L^{\infty }\left(
\mathbb{R}^{2}\right) \cap L^{1}\left( \mathbb{R}^{2}\right) $, in
particular, for the problem of evolution of vortex patches, where the
initial vorticity $q^{in}$ is assumed to be proportional to the
characteristic function of a bounded domain $\Omega ^{in}$, $%
q^{in}=q_{0}\chi _{\Omega ^{in}}$. Due to the conservation of the vorticity
along particle trajectories, the vorticity $q\left( t\right) $ remains a
characteristic function of an evolving in time domain $\Omega \left(
t\right) $. For the case where the boundary of the patch $\Omega ^{in}$
belongs to $C^{1,\gamma }$, $\gamma >0$, it was proved in \cite{a_C93} (see
also \cite{a_BC93}) that the Euler system has a unique solution $v\in
L_{loc}^{\infty }\left( \mathbb{R},Lip\left( \mathbb{R}^{2}\right) \right) $%
, and $\Omega \left( t\right) $ remains a bounded $C^{1,\gamma }$ domain. It
was proved in \cite{a_CW95} that the $L^{2}$ norm of the difference between
the solutions of NSE and the corresponding solution of the Euler system for
such initial data converges to zero, as the kinematic viscosity $\nu
\rightarrow 0$ (even though none of the solutions are in $L^{2}$). The rate
of convergence was improved to $\left( \nu t\right) ^{3/4}$ in \cite{a_AD04}%
, due to the fact that the vorticity of the vortex patch with $C^{1,\gamma }$
boundary is in fact in a Besov space $\dot{B}_{2,\infty }^{1/2}$, see also
\cite{a_M07} for a simpler proof and an extension to $\mathbb{R}^{3}$.

It has been shown in \cite{a_OS01} that there exists a unique global weak
solution to the Euler-$\alpha $ equations for initial vorticity in $\mathcal{%
M}({\mathbb{R}}^{2})$, the space of finite Radon measures on ${\mathbb{R}}%
^{2}$, with a unique Lagrangian flow map describing the evolution of
particles. In \cite{a_BLT09} we show the convergence, as $\alpha \rightarrow
0$, of the weak solutions of Euler-$\alpha $ equations with a distinguished
sign initial vorticity in $\mathcal{M}({\mathbb{R}}^{2})\cap
H_{loc}^{-1}\left( \mathbb{R}^{2}\right) $ (vortex sheet data) to those of
the 2D Euler equations.

Since the solution $q^{\alpha }$ of the Euler-$\alpha $ equation is
transported by the smoothed vector field $u^{\alpha }$, then in the Euler-$%
\alpha $ case the vortex patch is also transported by the flow. The global
existence and uniqueness results for the smooth vortex patch evolution under
the Euler-$\alpha $ equation can be obtained using arguments similar to
those presented in \cite{a_BLT09} and \cite[Chapter 8]{b_MB02}.
Specifically, if $q^{in}\left( x\right) =q_{0}\chi _{\Omega ^{in}}\left(
x\right) $ is a multiple of the characteristic function of a simply
connected bounded domain $\Omega ^{in}$, and the boundary $\partial \Omega
^{in}$ is in either one of the following spaces: $\mathrm{Lip}$, or $%
C^{1,\beta }$, $0\leq \beta \leq 1$, or $C^{2,\beta }$, $0\leq \beta <1$, or
$C^{n,\beta }$, $n\geq 3$, $0<\beta <1$, then the boundary of the vortex
patch remains in the same space as $\partial \Omega ^{in}$ for all times. We
describe this result in details in Appendix, Section \ref{sec:GlobalReg}.

In studying the convergence rate we use the following results.

\begin{proposition}
\label{prop_yudovich}Let $q^{in}\in L^{\infty }\left( \mathbb{R}^{2}\right)
\cap L^{1}\left( \mathbb{R}^{2}\right) $. Then there exist unique global
solutions $q$ and $q^{\alpha }$ of Euler and Euler-$\alpha $ equations %
\eqref{grp:EulerEqVortForm} and \eqref{grp:EulerAlphaEqVortForm}
respectively. Moreover, the $L^{p}$ norms of $q$ (of %
\eqref{grp:EulerEqVortForm}) and $q^{\alpha }$ (of %
\eqref{grp:EulerAlphaEqVortForm}) are conserved, namely, $\left\Vert q\left(
\cdot ,t\right) \right\Vert _{L^{p}}=\left\Vert q^{\alpha }\left( \cdot
,t\right) \right\Vert _{L^{p}}=\left\Vert q^{in}\right\Vert _{L^{p}},1\leq
p\leq \infty $. In addition, the velocities are bounded uniformly
\begin{equation}
\left\Vert v\left( \cdot ,t\right) \right\Vert _{L^{\infty }},\left\Vert
v^{\alpha }\left( \cdot ,t\right) \right\Vert _{L^{\infty }}\leq \left(
\left\Vert q^{in}\right\Vert _{L^{1}}\left\Vert q^{in}\right\Vert
_{L^{\infty }}\right) ^{1/2}.  \label{eq:vel_Linf_bound}
\end{equation}
\end{proposition}

The bounds on the velocities $v$ (of \eqref{grp:EulerEqVortForm}) and $%
v^{\alpha }$ (of \eqref{grp:EulerAlphaEqVortForm}) are a direct consequence
of the Biot-Savart law and the conservation of the $L^{1}$ and $L^{\infty }$
norms of vorticity.

\begin{proposition}
\label{prop_BC93} Let $q^{in}\left( x\right) =q_{0}\chi _{\Omega
^{in}}\left( x\right) $ be a multiple of the characteristic function of a
simply connected bounded domain $\Omega ^{in}$ with $C^{1,\gamma }$, $\gamma
\in \left( 0,1\right) $, boundary. Then the global solutions $v$ and $%
v^{\alpha }$ of Euler and Euler-$\alpha $ equations %
\eqref{grp:EulerEqVortForm} and \eqref{grp:EulerAlphaEqVortForm},
respectively, are in $L_{loc}^{\infty }\left( \mathbb{R};Lip\left( \mathbb{R}%
^{2}\right) \right) $ and for all $\alpha >0$, $t\in \mathbb{R}$
\begin{align}
& \left\Vert \nabla v\left( \cdot ,t\right) \right\Vert _{L^{\infty }}\leq
\left\Vert \nabla v^{in}\right\Vert _{L^{\infty }}e^{C\left\vert
t\right\vert },  \label{eq:grad_v_L_inf_bnd} \\
& \left\Vert \nabla v^{\alpha }\left( \cdot ,t\right) \right\Vert
_{L^{\infty }}\leq \left\Vert \nabla v^{in}\right\Vert _{L^{\infty
}}e^{C\left\vert t\right\vert },  \notag
\end{align}%
where $C=C\left( q^{in}\right) $, independent of $\alpha $, and the boundary
of the vortex patch remains $C^{1,\gamma }$ for all time.
\end{proposition}

For the Euler equations this result has been shown in \cite{a_BC93}. We
adopt their proof to show the uniform, in $\alpha $, bound on $\left\Vert
\nabla v^{\alpha }\left( \cdot ,t\right) \right\Vert _{L^{\infty }}$. The
vortex-patch problem is reformulated in terms of a scalar function $\varphi
^{\alpha }\left( x,t\right) $ that defines the patch boundary by $\Omega
^{\alpha }\left( t\right) =\left\{ \times \in \mathbb{R}^{2}|\varphi
^{\alpha }\left( x,t\right) >0\right\} $ and is convected with the flow by%
\begin{align*}
& \frac{\partial \varphi ^{\alpha }}{\partial t}+\left( {u}^{\alpha }\cdot
\nabla \right) \varphi ^{\alpha }=0, \\
& \varphi ^{\alpha }(x,0)=\varphi ^{in}(x).
\end{align*}%
To apply the method of the proof used in \cite{a_BC93}, the only ingredient
we need is to show that $\nabla {u}^{\alpha }$ is uniformly, in $\alpha $,
continuous in the tangential direction of the boundary. First, we recall
some properties of the kernel $K^{\alpha }$
\begin{equation}
K^{\alpha }\left( x\right) =\nabla ^{\perp }\Psi ^{\alpha }\left( \left\vert
x\right\vert \right) =\frac{x^{\perp }}{\left\vert x\right\vert }D\Psi
^{\alpha }\left( \left\vert x\right\vert \right) ,  \label{eq:K_alpha}
\end{equation}%
where
\begin{align}
& \Psi ^{\alpha }\left( r\right) =\frac{1}{2\pi }\left[ K_{0}\left( \frac{r}{%
\alpha }\right) +\log r\right] ,  \label{eq:DPsi} \\
& D\Psi ^{\alpha }(r)=\frac{d\Psi ^{\alpha }}{dr}(r)=\frac{1}{2\pi }\left[ -%
\frac{1}{\alpha }K_{1}\left( \frac{r}{\alpha }\right) +\frac{1}{r}\right] ,
\notag \\
& D^{2}\Psi ^{\alpha }\left( \left\vert x\right\vert \right) =\frac{1}{2\pi }%
\left[ \frac{1}{\alpha \left\vert x\right\vert }K_{1}\left( \frac{\left\vert
x\right\vert }{\alpha }\right) +\frac{1}{\alpha ^{2}}K_{0}\left( \frac{%
\left\vert x\right\vert }{\alpha }\right) -\frac{1}{\left\vert x\right\vert
^{2}}\right] ,  \notag \\
& D^{3}\Psi ^{\alpha }\left( r\right) =\frac{1}{2\pi }\left[ -\frac{2}{%
\alpha r^{2}}K_{1}\left( \frac{r}{\alpha }\right) -\frac{1}{\alpha ^{2}r}%
K_{0}\left( \frac{r}{\alpha }\right) -\frac{1}{\alpha ^{3}}K_{1}\left( \frac{%
r}{\alpha }\right) +\frac{2}{r^{3}}\right] .  \notag
\end{align}%
The functions $K_{0}$ and $K_{1}$ denote the modified Bessel functions of
the second kind of orders zero and one, respectively. For details on Bessel
functions, see, e.g., \cite{b_W44}. Derivatives of $\Psi ^{\alpha }$ decay
to zero as $\frac{r}{\alpha }\rightarrow \infty $; and as $\frac{r}{\alpha }%
\rightarrow 0$ satisfy%
\begin{align}
& D\Psi ^{\alpha }\left( r\right) =-\frac{1}{4\pi }\frac{r}{\alpha ^{2}}\log
\frac{r}{\alpha }+O\left( \frac{r}{\alpha ^{2}}\right) ,
\label{eq:Der_psi_at_origin} \\
& D^{2}\Psi ^{\alpha }\left( r\right) =-\frac{1}{4\pi }\frac{1}{\alpha ^{2}}%
\log \frac{r}{\alpha }+O\left( \frac{1}{\alpha ^{2}}\right) ,  \notag \\
& D^{3}\Psi ^{\alpha }\left( r\right) =-\frac{1}{4\pi }\frac{1}{r\alpha ^{2}}%
+O\left( \frac{r}{\alpha ^{4}}\log \frac{r}{\alpha }\right) ,  \notag
\end{align}%
where the constants in the big O are independent of $\alpha $. The filtered
velocity gradient is given by%
\begin{equation*}
\nabla {u}^{\alpha }\left( x,t\right) =\int_{\mathbb{R}^{2}}\nabla K^{\alpha
}\left( x-y\right) q^{\alpha }\left( y,t\right) dy
\end{equation*}%
As in \cite{a_BC93}, denote by $W$ a divergence free vector field which is
tangent to $\partial \Omega $, $W=\nabla ^{\perp }\varphi $. Then
\begin{equation*}
\nabla {u}^{\alpha }\left( x\right) W\left( x\right) =\int_{\mathbb{R}%
^{2}}\nabla K^{\alpha }\left( x-y\right) q^{\alpha }\left( y\right) \left[
W\left( x\right) -W\left( y\right) \right] dy,
\end{equation*}%
and we have the following result corresponding to Corollary 1 and Lemma in
the appendix of \cite{a_BC93}.

\begin{lemma}
For $\gamma \in \left( 0,1\right) $%
\begin{equation*}
\left\Vert \nabla {u}^{\alpha }W\right\Vert _{C^{0,\gamma }}\leq C\left\Vert
\nabla {v}^{\alpha }\right\Vert _{L^{\infty }}\left\Vert W\right\Vert
_{C^{0,\gamma }}.
\end{equation*}
\end{lemma}

\begin{proof}
We stress that all the constants $C$ are independent of $\alpha $. We write%
\begin{align*}
\nabla u^{\alpha }\left( x\right) W\left( x\right) & -\nabla u^{\alpha
}\left( x+h\right) W\left( x+h\right) \\
& =\int_{\mathbb{R}^{2}}\nabla K^{\alpha }\left( x-y\right) q^{\alpha
}\left( y\right) \left[ W\left( x\right) -W\left( y\right) \right] dy \\
& -\int_{\mathbb{R}^{2}}\nabla K^{\alpha }\left( x+h-y\right) q^{\alpha
}\left( y\right) \left[ W\left( x+h\right) -W\left( y\right) \right] dy \\
& =\int_{\left\vert x-y\right\vert <2\left\vert h\right\vert }\nabla
K^{\alpha }\left( x-y\right) q^{\alpha }\left( y\right) \left[ W\left(
x\right) -W\left( y\right) \right] dy \\
& -\int_{\left\vert x-y\right\vert <2\left\vert h\right\vert }\nabla
K^{\alpha }\left( x+h-y\right) q^{\alpha }\left( y\right) \left[ W\left(
x+h\right) -W\left( y\right) \right] dy \\
& +\int_{\left\vert x-y\right\vert \geq 2\left\vert h\right\vert }\nabla
K^{\alpha }\left( x-y\right) q^{\alpha }\left( y\right) \left[ W\left(
x\right) -W\left( x+h\right) \right] dy \\
& +\int_{\left\vert x-y\right\vert \geq 2\left\vert h\right\vert }\left[
\nabla K^{\alpha }\left( x-y\right) -\nabla K^{\alpha }\left( x+h-y\right) %
\right] q^{\alpha }\left( y\right) \left[ W\left( x+h\right) -W\left(
y\right) \right] dy \\
& =I_{1}+I_{2}+I_{3}+I_{4}.
\end{align*}%
Using the fact that $\left\vert \nabla K^{\alpha }\left( x\right)
\right\vert \leq \frac{C}{\left\vert x\right\vert ^{2}}$ (cf. %
\eqref{eq:K_alpha}, \eqref{eq:DPsi}), we obtain that $\left\vert
I_{1}\right\vert ,\left\vert I_{2}\right\vert \leq C\left\Vert q^{\alpha
}\right\Vert _{L^{\infty }}\left\Vert W\right\Vert _{C^{0,\gamma }}h^{\gamma
}$, also due to $\left\vert D^{2}K^{\alpha }\left( x\right) \right\vert \leq
\frac{C}{\left\vert x\right\vert ^{3}}$ (cf. \eqref{eq:K_alpha}, %
\eqref{eq:DPsi}), we have $\left\vert I_{4}\right\vert \leq C\left( \gamma
\right) \left\Vert q^{\alpha }\right\Vert _{L^{\infty }}\left\Vert
W\right\Vert _{C^{0,\gamma }}h^{\gamma }$. To bound the term $I_{3}$ we
consider two cases $\left\vert h\right\vert \leq \alpha $ and $\left\vert
h\right\vert >\alpha $ separately. We have%
\begin{equation*}
\left\vert I_{3}\right\vert \leq \left\Vert W\right\Vert _{C^{0,\gamma
}}h^{\gamma }\left\vert J\right\vert ,
\end{equation*}%
where $J=\int_{\left\vert x-y\right\vert \geq 2\left\vert h\right\vert
}\nabla K^{\alpha }\left( x-y\right) q^{\alpha }\left( y\right) dy$. First,
let $\left\vert h\right\vert \leq \alpha $, write $J$ as%
\begin{align*}
J& =\int_{\mathbb{R}^{2}}\nabla K^{\alpha }\left( x-y\right) q^{\alpha
}\left( y\right) dy-\int_{\left\vert x-y\right\vert <2\left\vert
h\right\vert }\nabla K^{\alpha }\left( x-y\right) q^{\alpha }\left( y\right)
dy \\
& =\nabla u^{\alpha }\left( x\right) -J_{1}.
\end{align*}%
By \eqref{eq:K_alpha}-\eqref{eq:Der_psi_at_origin}, we obtain%
\begin{align*}
\left\vert J_{1}\right\vert & \leq \left\Vert q^{\alpha }\right\Vert
_{L^{\infty }}\int_{\left\vert x-y\right\vert <2\left\vert h\right\vert
}\left\vert \nabla K^{\alpha }\left( x-y\right) \right\vert dy \\
& \leq \left\Vert q^{\alpha }\right\Vert _{L^{\infty }}\int_{\frac{%
\left\vert x-y\right\vert }{\alpha }<2}\left( \frac{1}{4\pi }\frac{1}{\alpha
^{2}}\left\vert \log \frac{\left\vert x-y\right\vert }{\alpha }\right\vert +%
\frac{C}{\alpha ^{2}}\right) dy \\
& \leq C\left\Vert q^{\alpha }\right\Vert _{L^{\infty }}\text{.}
\end{align*}%
Now, let $\left\vert h\right\vert >\alpha $, then we write%
\begin{align*}
J& =\int_{\left\vert x-y\right\vert \geq 2\left\vert h\right\vert }\left(
\nabla K^{\alpha }\left( x-y\right) -\nabla K\left( x-y\right) \right)
q^{\alpha }\left( y\right) dy+\int_{\left\vert x-y\right\vert \geq
2\left\vert h\right\vert }\nabla K\left( x-y\right) q^{\alpha }\left(
y\right) dy \\
& =J_{2}+J_{3}.
\end{align*}%
A bound for $\left\vert J_{3}\right\vert $ is obtained by using a lemma due
to Cotlar (see \cite{b_T86} p. 291)%
\begin{align*}
\left\vert J_{3}\right\vert & \leq C\left( \left\Vert \mathrm{p.v.}\int_{%
\mathbb{R}^{2}}\nabla K\left( x-y\right) q^{\alpha }\left( y\right)
dy\right\Vert _{L^{\infty }}+\left\Vert q^{\alpha }\right\Vert _{L^{\infty
}}\right) \\
& =C\left( \left\Vert \nabla v^{\alpha }\right\Vert _{L^{\infty
}}+\left\Vert q^{\alpha }\right\Vert _{L^{\infty }}\right) .
\end{align*}%
We bound $\left\vert J_{2}\right\vert $ using the facts that $K_{0}$,$%
K_{1}\geq 0$ satisfy $\int_{\mathbb{R}^{2}}\frac{1}{\alpha ^{2}}K_{0}\left(
\frac{\left\vert x\right\vert }{\alpha }\right) dx=2\pi ,$ $\int_{\left\vert
x\right\vert >2\alpha }\frac{1}{\alpha \left\vert x\right\vert }K_{1}\left(
\frac{\left\vert x\right\vert }{\alpha }\right) dx=\pi K_{0}\left( 2\right) $%
, and we obtain
\begin{align*}
\left\vert J_{2}\right\vert & \leq C\left\Vert q^{\alpha }\right\Vert
_{L^{\infty }}\int_{\left\vert x-y\right\vert >2\alpha }\left[ \frac{1}{%
\alpha \left\vert x-y\right\vert }K_{1}\left( \frac{\left\vert
x-y\right\vert }{\alpha }\right) +\frac{1}{\alpha ^{2}}K_{0}\left( \frac{%
\left\vert x-y\right\vert }{\alpha }\right) \right] dy \\
& \leq C\left\Vert q^{\alpha }\right\Vert _{L^{\infty }}.
\end{align*}%
To conclude, in both cases we have
\begin{equation*}
\left\vert I_{3}\right\vert \leq C\left\Vert W\right\Vert _{C^{0,\gamma
}}h^{\gamma }\left( \left\Vert \nabla v^{\alpha }\right\Vert _{L^{\infty
}}+\left\Vert q^{\alpha }\right\Vert _{L^{\infty }}\right) .
\end{equation*}
\end{proof}

Next we briefly recall the definition of homogeneous Besov spaces. Let $%
\mathcal{S}\left( \mathbb{R}^{n}\right) $ be the Schwartz space and denote
by $\mathcal{Z}^{\prime }\left( \mathbb{R}^{n}\right) $ the dual space of $%
\mathcal{Z}\left( \mathbb{R}^{n}\right) =\left\{ f\in S\left( \mathbb{R}%
^{n}\right) :D^{\beta }\hat{f}\left( 0\right) =0\text{ for every multi-index
}\beta \in \mathbb{N}^{n}\right\} $, it can also be identified as the
quotient space of $\mathcal{S}^{\prime }/\mathcal{P}$, where $\mathcal{P}$
is the collection of all polynomials. Here $\hat{\varphi}$ and $\mathcal{F}%
\left( \varphi \right) $ denote the Fourier transform of $\varphi $ in $%
\mathbb{R}^{n}$.

We recall the Littlewood-Paley decomposition. Choose a radial function $%
\varphi \in \mathcal{S}$ such that $\hat{\varphi}\in C_{0}^{\infty }\left(
\mathbb{R}^{n}\backslash \left\{ 0\right\} \right) $ satisfies $\supp\hat{%
\varphi}\subset \left\{ \frac{3}{4}\leq \left\vert \xi \right\vert \leq
\frac{8}{3}\right\} $ and $\sum_{q=-\infty }^{\infty }\hat{\varphi}%
_{q}\left( \xi \right) =1$ for $\xi \neq 0$, where $\hat{\varphi}_{q}$%
\thinspace $\left( \xi \right) =\hat{\varphi}\left( 2^{-q}\xi \right) $,
that is, $\varphi _{q}$\thinspace $\left( x\right) =2^{qn}\varphi \left(
2^{q}x\right) $. For $q\in \mathbb{Z}$ one defines the dyadic blocks by $%
\Delta _{q}f=\mathcal{F}^{-1}\left( \mathcal{F}\left( \varphi _{q}\right)
\mathcal{F}\left( f\right) \right) =\varphi _{q}\ast f$. The formal
decomposition $f=\sum_{q=-\infty }^{\infty }\Delta _{q}f$ holds true modulo
polynomials and is called the homogeneous Littlewood-Paley decomposition.
For $s\in \mathbb{R}$, $1\leq r\leq \infty $, the homogeneous Besov space is
defined as
\begin{equation*}
\dot{B}_{2,r}^{s}=\left\{ f\in \mathcal{Z}^{\prime }\left( \mathbb{R}%
^{n}\right) :|\left\Vert f\right\Vert _{\dot{B}_{2,r}^{s}}<\infty \right\} ,
\end{equation*}%
where $\left\Vert f\right\Vert _{\dot{B}_{2,r}^{s}}=\left( \sum_{q=-\infty
}^{\infty }2^{rsq}\left\Vert \Delta _{q}f\right\Vert _{L^{2}}^{r}\right)
^{1/r}$ for $r\in \left[ 1,\infty \right) $ and $\left\Vert f\right\Vert _{%
\dot{B}_{2,\infty }^{s}}=\sup_{q\in \mathbb{Z}}\left( 2^{sq}\left\Vert
\Delta _{q}f\right\Vert _{L^{2}}\right) $. For the detailed description of
Besov space see, e.g., \cite{b_BL76,b_T83,b_P76,b_L02}.

We use the following results:

\begin{lemma}
\label{lem:interpolation_ineq}\cite{a_M07} For $0<\beta <1$
\begin{equation}
\left\Vert f\right\Vert _{\dot{B}_{2,1}^{1-\beta }}\leq C\left\Vert
f\right\Vert _{L^{2}}^{\beta }\left\Vert \nabla f\right\Vert
_{L^{2}}^{1-\beta }.  \label{eq:interpolation_ineq}
\end{equation}
\end{lemma}

\begin{lemma}
\label{lem:BesovSpaceDerIneq} Let $s\in \mathbb{R}$, $\beta \geq 0$. There
exist constants $c$ and $C$ such that%
\begin{equation}
c^{-\beta }\left\Vert f\right\Vert _{B_{2,\infty }^{s+\beta }}\leq
\left\Vert \left( -\Delta \right) ^{\frac{\beta }{2}}f\right\Vert
_{B_{2,\infty }^{s}}\leq C^{\beta }\left\Vert f\right\Vert _{B_{2,\infty
}^{s+\beta }}  \label{eq:Besov_D}
\end{equation}
\end{lemma}

\begin{proof}
This lemma follows directly from definition of the Besov spaces since
\begin{equation*}
\left\Vert \Delta _{q}\left( \left( -\Delta \right) ^{\frac{\beta }{2}%
}f\right) \right\Vert _{L^{2}}=\left\Vert \mathcal{F}\left( \varphi
_{q}\right) \mathcal{F}\left( \left( -\Delta \right) ^{\frac{\beta }{2}%
}f\right) \right\Vert _{L^{2}}=\left\Vert \hat{\varphi}\left( 2^{-q}\xi
\right) \left\vert \xi \right\vert ^{\beta }\hat{f}\left( \xi \right)
\right\Vert _{L^{2}}
\end{equation*}%
and because $\hat{\varphi}\left( 2^{-q}\xi \right) $ is supported in $%
\left\{ \xi \in \mathbb{R}^{n}:\frac{3}{4}2^{q}\leq \left\vert \xi
\right\vert \leq \frac{8}{3}2^{q}\right\} $.
\end{proof}

\begin{lemma}
Let $0\leq \beta \leq 2$, $g\in \dot{B}_{2,r}^{\beta }$ and $f=\left(
1-\alpha ^{2}\Delta \right) ^{-1}g$. Then%
\begin{equation}
\alpha ^{\beta }\left\Vert \left( -\Delta \right) ^{\beta /2}f\right\Vert _{%
\dot{B}_{2,r}^{\beta }}\leq \left\Vert g\right\Vert _{\dot{B}_{2,r}^{\beta }}
\label{eq:alpha_Besov_grad_u}
\end{equation}
\end{lemma}

\begin{proof}
We have%
\begin{align*}
\alpha ^{\beta }\left\Vert \Delta _{q}\left( \left( -\Delta \right) ^{\beta
/2}f\right) \right\Vert _{L^{2}}& =\alpha ^{\beta }\left\Vert \mathcal{F}%
\left( \varphi _{q}\right) \mathcal{F}\left( \left( -\Delta \right) ^{\frac{%
\beta }{2}}f\right) \right\Vert _{L^{2}} \\
& =\alpha ^{\beta }\left\Vert \mathcal{F}\left( \varphi _{q}\right) \frac{%
\left\vert \xi \right\vert ^{\beta }\hat{g}\left( \xi \right) }{1+\alpha
^{2}\left\vert \xi \right\vert ^{2}}\right\Vert _{L^{2}} \\
& =\left( \int_{\mathbb{R}^{n}}\frac{\left( \alpha ^{2}\left\vert \xi
\right\vert ^{2}\right) ^{\beta }}{\left( 1+\alpha ^{2}\left\vert \xi
\right\vert ^{2}\right) ^{2}}\left( \hat{\varphi}\left( 2^{-q}\xi \right)
\hat{g}\left( \xi \right) \right) ^{2}d\xi \right) ^{1/2} \\
& \leq \left( \left( \sup_{y\geq 0}\frac{y^{\beta }}{\left( 1+y\right) ^{2}}%
\right) \int_{\mathbb{R}^{n}}\left( \hat{\varphi}\left( 2^{-q}\xi \right)
\hat{g}\left( \xi \right) \right) ^{2}d\xi \right) ^{1/2} \\
& \leq \left\Vert \Delta _{q}g\right\Vert _{L^{2}}.
\end{align*}%
From which \eqref{eq:alpha_Besov_grad_u} follows.
\end{proof}

The two following propositions show that the vortex patch vorticity, with $%
C^{1,\gamma }$ boundary, evolving under Euler-$\alpha $ equations is in a
homogeneous Besov space $\dot{B}_{2,\infty }^{1/2}$.

\begin{proposition}
\label{prop:vort_is_besov}\cite{a_M07} If $\Omega $ is a $C^{1,\gamma }$, $%
\gamma >0$, bounded domain, then $\chi _{\Omega }\in \dot{B}_{2,\infty
}^{1/2}$.
\end{proposition}

\begin{proposition}
\label{prop:vort_is_besov_in_time}Let $v\in L^{\infty }\left( \left[ 0,T%
\right] ;Lip\left( \mathbb{R}^{2}\right) \right) $ with $\diver v=0$. Let $%
\varphi $ solve%
\begin{align*}
& \frac{\partial \varphi }{\partial t}+\left( {v}\cdot \nabla \right)
\varphi =0, \\
& \varphi (x,0)=\varphi ^{in}(x),
\end{align*}%
with $\varphi ^{in}\in \dot{B}_{2,\infty }^{1/2}$. Then%
\begin{equation*}
\left\Vert \varphi \left( \cdot ,t\right) \right\Vert _{\dot{B}_{2,\infty
}^{1/2}}\leq \left\Vert \varphi ^{in}\right\Vert _{\dot{B}_{2,\infty
}^{1/2}}e^{C\int_{0}^{t}\left\Vert \nabla v\left( \cdot ,\tau \right)
\right\Vert _{L^{\infty }}d\tau }.
\end{equation*}
\end{proposition}

This result is a straightforward adaptation of the Proposition 3.1 of \cite%
{a_AD04} to the homogeneous Besov spaces.

Propositions \ref{prop_BC93}, \ref{prop:vort_is_besov} and \ref%
{prop:vort_is_besov_in_time} imply the following result:

\begin{corollary}
\label{prop:vort_Besov} Let $q^{in}\left( x\right) =q_{0}\chi _{\Omega
^{in}}\left( x\right) $ be a multiple of the characteristic function of a
simply connected bounded domain $\Omega ^{in}$ with $C^{1,\gamma }$, $\gamma
>0$, boundary. Then the global solutions $q$ and $q^{\alpha }$ of Euler and
Euler-$\alpha $ equations \eqref{grp:EulerEqVortForm} and %
\eqref{grp:EulerAlphaEqVortForm}, respectively, satisfy for all $\alpha >0$,
$t\in \mathbb{R}^{+}$
\begin{eqnarray*}
\left\Vert q\left( \cdot ,t\right) \right\Vert _{\dot{B}_{2,\infty }^{1/2}}
&\leq &\left\Vert q^{in}\right\Vert _{\dot{B}_{2,\infty
}^{1/2}}e^{C\int_{0}^{t}\left\Vert \nabla v\left( \cdot ,\tau \right)
\right\Vert _{L^{\infty }}d\tau }, \\
\left\Vert q^{\alpha }\left( \cdot ,t\right) \right\Vert _{\dot{B}_{2,\infty
}^{1/2}} &\leq &\left\Vert q^{in}\right\Vert _{\dot{B}_{2,\infty
}^{1/2}}e^{C\int_{0}^{t}\left\Vert \nabla u^{\alpha }\left( \cdot ,\tau
\right) \right\Vert _{L^{\infty }}d\tau }.
\end{eqnarray*}
\end{corollary}

Next we state a lemma that play an important role in estimating the
convergence rate.

\begin{lemma}
Let $q^{in}$ $=q_{0}\chi _{\Omega ^{in}}$, with $\Omega ^{in}$ being a
simply connected bounded domain $\Omega ^{in}$ with $C^{1,\gamma }$
boundary, $\gamma \in \left( 0,1\right) $. Then for all $T\in \left[
0,\infty \right) $ the solution $u^{\alpha }$ of the Euler-$\alpha $
equations \eqref{grp:EulerAlphaEq} satisfies
\begin{equation}
\alpha ^{1/2}\left\Vert \Delta u^{\alpha }\left( \cdot ,t\right) \right\Vert
_{L^{2}}\leq C\left\Vert q^{\alpha }\right\Vert _{L^{\infty }\left( \left[
0,T\right] ,\dot{B}_{2,\infty }^{1/2}\right) }
\label{eq:alpha_square_laplacian_u_bound}
\end{equation}%
for all $t\in \left[ 0,T\right] $.
\end{lemma}

\begin{proof}
Let $0<\beta <\frac{1}{2}$. Since $\left( \dot{B}_{2,1}^{-\beta }\right)
^{\prime }=\dot{B}_{2,\infty }^{\beta }$, for $\beta \in \mathbb{R}$ (see,
e.g., \cite{b_T83}), we have%
\begin{equation*}
\alpha ^{\beta +1/2}\left\Vert \Delta u^{\alpha }\right\Vert
_{L^{2}}^{2}=\alpha ^{\beta +1/2}\int_{\mathbb{R}^{n}}\Delta u^{\alpha
}\left( x\right) \cdot \Delta u^{\alpha }\left( x\right) dx\leq C\alpha
^{\beta +1/2}\left\Vert \Delta u^{\alpha }\right\Vert _{\dot{B}_{2,\infty
}^{\beta }}\left\Vert \Delta u^{\alpha }\right\Vert _{\dot{B}_{2,1}^{-\beta
}}.
\end{equation*}%
We stress that all the constants $C$ are independent of $\alpha $ and $\beta
$. Now, by \eqref{eq:alpha_Besov_grad_u} and \eqref{eq:Besov_D}%
\begin{equation*}
\alpha ^{\beta +1/2}\left\Vert \Delta u^{\alpha }\right\Vert _{\dot{B}%
_{2,\infty }^{\beta }}
\end{equation*}
\begin{equation*}
\alpha ^{\beta +1/2}\left\Vert \Delta u^{\alpha }\right\Vert _{\dot{B}%
_{2,\infty }^{\beta }}\leq \left\Vert \left( -\Delta \right) ^{3/4-\beta
/2}v^{\alpha }\right\Vert _{\dot{B}_{2,\infty }^{\beta }}\leq C\left\Vert
\left( -\Delta \right) ^{1/2}v^{\alpha }\right\Vert _{\dot{B}_{2,\infty
}^{1/2}}\leq C\left\Vert q^{\alpha }\right\Vert _{\dot{B}_{2,\infty }^{1/2}},
\end{equation*}%
and by \eqref{eq:interpolation_ineq}%
\begin{align*}
\left\Vert \Delta u^{\alpha }\right\Vert _{\dot{B}_{2,1}^{-\beta }}& \leq
C\left\Vert \nabla u^{\alpha }\right\Vert _{L^{2}}^{\beta }\left\Vert \Delta
u^{\alpha }\right\Vert _{L^{2}}^{1-\beta } \\
& \leq C\left\Vert q^{in}\right\Vert _{L^{2}}^{\beta }\left\Vert \Delta
u^{\alpha }\right\Vert _{L^{2}}^{1-\beta }.
\end{align*}%
Summing up we have%
\begin{equation*}
\alpha ^{\beta +1/2}\left\Vert \Delta u^{\alpha }\right\Vert
_{L^{2}}^{2}\leq C\left\Vert q^{\alpha }\right\Vert _{\dot{B}_{2,\infty
}^{1/2}}\left\Vert q^{in}\right\Vert _{L^{2}}^{\beta }\left\Vert \Delta
u^{\alpha }\right\Vert _{L^{2}}^{1-\beta },
\end{equation*}%
hence
\begin{equation*}
\alpha ^{\frac{\beta +1/2}{1+\beta }}\left\Vert \Delta u^{\alpha
}\right\Vert _{L^{2}}\leq C\left\Vert q^{\alpha }\right\Vert _{\dot{B}%
_{2,\infty }^{1/2}}^{1/\left( 1+\beta \right) }\left\Vert q^{in}\right\Vert
_{L^{2}}^{\beta /\left( 1+\beta \right) },
\end{equation*}%
and taking the limit as $\beta \rightarrow 0$ we obtain%
\begin{equation*}
\alpha ^{\frac{1}{2}}\left\Vert \Delta u^{\alpha }\right\Vert _{L^{2}}\leq
C\left\Vert q^{\alpha }\right\Vert _{\dot{B}_{2,\infty }^{1/2}}.
\end{equation*}
\end{proof}

In the next theorem we show that the solution of the Euler-$\alpha $
equations differs from the solution of the Euler equations by order $\left(
\alpha ^{2}\right) ^{3/4}$, both having the same vortex patch as initial
data.

\begin{theorem}
Let $q^{in}=q_{0}\chi _{\Omega ^{in}}$, where $\chi _{\Omega ^{in}}$ is a
characteristic function of a simply connected bounded domain $\Omega ^{in}$
with $C^{1,\gamma }$, $\gamma \in \left( 0,1\right) $, boundary, and let $%
v,v^{\alpha }$ be the solutions of the Euler and the Euler-$\alpha $
equations \eqref{grp:EulerEq} and \eqref{grp:EulerAlphaEq}, respectively,
with initial data $q^{in}$. Then the difference $w^{\alpha }=v-v^{\alpha }$
is square-integrable and obeys the estimate%
\begin{equation*}
\left\Vert w^{\alpha }\left( \cdot ,t\right) \right\Vert _{L^{2}}^{2}\leq
\alpha ^{3/2}CV\left( t\right) \left\Vert q^{in}\right\Vert _{\dot{B}%
_{2,\infty }^{1/2}}e^{CV\left( t\right) },
\end{equation*}%
where%
\begin{equation*}
V\left( t\right) =\int_{0}^{t}\left( \left\Vert \nabla v\left( \cdot ,\tau
\right) \right\Vert _{L^{\infty }}+\left\Vert \nabla u^{\alpha }\left( \cdot
,\tau \right) \right\Vert _{L^{\infty }}\right) d\tau .
\end{equation*}%
In particular, there exists a constant $C=C\left( q^{in}\right) $ such that
\begin{equation*}
\left\Vert w^{\alpha }\left( \cdot ,t\right) \right\Vert _{L^{2}}\leq \alpha
^{3/2}Ce^{Ce^{Ct}}.
\end{equation*}
\end{theorem}

\begin{proof}
The difference $w^{\alpha }=v-v^{\alpha }$ satisfies the following equation%
\begin{align}
\frac{\partial w^{\alpha }}{\partial t}& +(w^{\alpha }\cdot \nabla
)v+v_{j}\nabla w_{j}^{\alpha }+({u}^{\alpha }\cdot \nabla )w^{\alpha
}+w_{_{j}}^{\alpha }\nabla {u}_{j}^{\alpha }
\label{eq:diff_eq_vortex_patch_case} \\
& -(\alpha ^{2}\Delta {u}^{\alpha }\cdot \nabla )v-v_{j}\nabla \left( \alpha
^{2}\Delta {u}_{j}^{\alpha }\right) -v_{j}\nabla v_{j}+\nabla \left(
p-p^{\alpha }\right) =0.  \notag
\end{align}%
We remark that $v$ and $v^{\alpha }$, obtained from the vortex patch
vorticity by convolution with the Biot-Savart kernel are not in $L^{2}$,
however, the difference $w^{\alpha }\in L^{\infty }\left( \left[ 0,T\right]
,L^{2}\left( \mathbb{R}^{2}\right) \right) $, since $\curl w^{\alpha }$ is
compactly supported and $\int_{\mathbb{R}^{2}}\curl w^{\alpha }\left(
x,t\right) dx=\int_{\mathbb{R}^{2}}\left( q^{in}-q^{in,\alpha }\right)
\left( x\right) dx=0$ (this could be seen by using an asymptotic expansion
of the kernel, see, e.g., \cite{b_MB02}, p.321). Also, all the terms in %
\eqref{eq:diff_eq_vortex_patch_case} are in $L^{\infty }\left( \left[ 0,T%
\right] ,L^{2}\left( \mathbb{R}^{2}\right) \right) $, since the velocities
are in $L^{\infty }\left( \left[ 0,T\right] ,L^{\infty }\left( \mathbb{R}%
^{2}\right) \right) $, and their gradients, as well as gradients of the
pressures, are in $L^{\infty }\left( \left[ 0,T\right] ,L^{2}\left( \mathbb{R%
}^{2}\right) \right) $.

We take an $\,L^{2}$-inner product of \eqref{eq:diff_eq_vortex_patch_case}
with $w^{\alpha }\left( x,t\right) $ and obtain that%
\begin{equation}
\frac{1}{2}\frac{d}{dt}\left\Vert w^{\alpha }\left( \cdot ,t\right)
\right\Vert _{L^{2}}^{2}\leq I_{1}+I_{2}+I_{3},
\label{eq:vortexPatch_diff_bnd}
\end{equation}%
where
\begin{align*}
I_{1}& =\left\vert \left( w_{_{j}}^{\alpha }\nabla {u}_{j}^{\alpha
},w^{\alpha }\right) \right\vert , \\
I_{2}& =\alpha ^{2}\left\vert \left( (\Delta {u}^{\alpha }\cdot \nabla
)v,w^{\alpha }\right) \right\vert +\alpha ^{2}\left\vert \left( (w^{\alpha
}\cdot \nabla )v,\Delta {u}^{\alpha }\right) \right\vert , \\
I_{3}& =\left\vert \left( \nabla \left( p-p^{\alpha }\right) ,w^{\alpha
}\right) \right\vert .
\end{align*}%
here we used the identity, for $h$ divergence free,%
\begin{equation*}
\left( (f\cdot \nabla )g+g_{j}\nabla f_{j},h\right) =\left( (f\cdot \nabla
)g,h\right) -\left( (h\cdot \nabla )g,f\right) .
\end{equation*}%
The first term is estimated by
\begin{equation*}
I_{1}\leq \left\Vert w^{\alpha }\right\Vert _{L^{2}}^{2}\left\Vert \nabla
u^{\alpha }\right\Vert _{L^{\infty }}.
\end{equation*}%
For the second term by \eqref{eq:alpha_square_laplacian_u_bound} we obtain
\begin{align*}
I_{2}& \leq \alpha ^{2}\left\Vert \Delta {u}^{\alpha }\right\Vert
_{L^{2}}\left\Vert \nabla v\right\Vert _{L^{\infty }}\left\Vert w^{\alpha
}\right\Vert _{L^{2}} \\
& \leq C\alpha ^{3/2}\left\Vert q^{\alpha }\right\Vert _{L^{\infty }\left(
\left[ 0,T\right] ,\dot{B}_{2,\infty }^{1/2}\right) }\left\Vert \nabla
v\right\Vert _{L^{\infty }}\left\Vert w^{\alpha }\right\Vert _{L^{2}}.
\end{align*}%
It remains to estimate the third term. We remark that the pressure is
determined uniquely up to a constant, and in order to ensure that $p\in L^{2}
$, we require a side condition $\int_{\mathbb{R}^{2}}p\left( x\right) dx=0$.
Notice that the mean free pressure of Euler equations satisfy
\begin{equation*}
p=R_{i}R_{j}\left( v_{i}v_{j}\right)
\end{equation*}%
where $R_{i}=\left( -\Delta \right) ^{-1/2}\frac{\partial }{\partial x_{i}}$
is the Riesz transform (see, e.g., \cite{b_G04}), hence by the properties of
the Riesz transform and the Sobolev embedding theorem we have%
\begin{equation*}
\left\Vert p\right\Vert _{L^{2}}=\left\Vert v_{i}v_{j}\right\Vert
_{L^{2}}\leq \left\Vert v\right\Vert _{L^{4}}^{2}\leq C\left\Vert
v\right\Vert _{W^{1,4/3}}^{2}\leq C\left\Vert q^{in}\right\Vert
_{L^{4/3}}^{2}.
\end{equation*}%
For the Euler-$\alpha $ equation we have
\begin{equation*}
-\Delta p^{\alpha }=\frac{\partial }{\partial x_{i}}\frac{\partial }{%
\partial x_{j}}\left( u_{i}^{\alpha }v_{j}^{\alpha }\right) +\frac{1}{2}%
\frac{\partial ^{2}}{\partial x_{i}^{2}}\left( u_{j}^{\alpha }\right) ^{2}-%
\frac{\partial }{\partial x_{i}}\left( \alpha ^{2}\Delta u_{j}^{\alpha }%
\frac{\partial }{\partial x_{i}}u_{j}^{\alpha }\right) ,
\end{equation*}%
we write%
\begin{equation*}
p^{\alpha }=p_{1}^{\alpha }+p_{2}^{\alpha },
\end{equation*}%
where the mean free pressure $p_{1}^{\alpha }$, $\int_{\mathbb{R}%
^{2}}p_{1}^{\alpha }\left( x\right) dx=0$, satisfies
\begin{equation*}
-\Delta p_{1}^{\alpha }=\frac{\partial }{\partial x_{i}}\frac{\partial }{%
\partial x_{j}}\left( u_{i}^{\alpha }v_{j}^{\alpha }\right) +\frac{1}{2}%
\frac{\partial ^{2}}{\partial x_{i}^{2}}\left( u_{j}^{\alpha }\right) ^{2},
\end{equation*}%
and $p_{2}^{\alpha }$ satisfies%
\begin{equation*}
-\Delta p_{2}^{\alpha }=-\alpha ^{2}\frac{\partial }{\partial x_{i}}\left(
\Delta u_{j}^{\alpha }\frac{\partial }{\partial x_{i}}u_{j}^{\alpha }\right)
,
\end{equation*}%
As described above for the pressure of the Euler equations, we have
\begin{equation*}
\left\Vert p_{1}^{\alpha }\right\Vert _{L^{2}}\leq C\left\Vert
q^{in}\right\Vert _{L^{4/3}}^{2}.
\end{equation*}%
For the $p_{2}^{\alpha }$ we have
\begin{equation*}
\left\Vert \nabla p_{2}^{\alpha }\right\Vert _{L^{2}}=\alpha ^{2}\left\Vert
\Delta u_{j}^{\alpha }\frac{\partial }{\partial x_{i}}u_{j}^{\alpha
}\right\Vert _{L^{2}}\leq \alpha ^{2}\left\Vert \Delta u^{\alpha
}\right\Vert _{L^{2}}\left\Vert \nabla u^{\alpha }\right\Vert _{L^{\infty }}.
\end{equation*}%
Hence, using also \eqref{eq:alpha_square_laplacian_u_bound}, we obtain%
\begin{equation*}
I_{3}\leq \left\Vert \nabla p_{2}^{\alpha }\right\Vert _{L^{2}}\left\Vert
w^{\alpha }\right\Vert _{L^{2}}\leq C\alpha ^{3/2}\left\Vert q^{\alpha
}\right\Vert _{L^{\infty }\left( \left[ 0,T\right] ,\dot{B}_{2,\infty
}^{1/2}\right) }\left\Vert \nabla u^{\alpha }\right\Vert _{L^{\infty
}}\left\Vert w^{\alpha }\right\Vert _{L^{2}}.
\end{equation*}%
From the above and \eqref{eq:vortexPatch_diff_bnd}, and by using Gr\"{o}%
nwall lemma, we obtain%
\begin{equation*}
\left\Vert w^{\alpha }\left( \cdot ,t\right) \right\Vert _{L^{2}}\leq
C\alpha ^{3/2}\left\Vert q^{in}\right\Vert _{\dot{B}_{2,\infty
}^{1/2}}V\left( t\right) e^{CV\left( t\right) },
\end{equation*}%
where $V\left( t\right) =\int_{0}^{t}\left( \left\Vert \nabla v\left( \cdot
,\tau \right) \right\Vert _{L^{\infty }}+\left\Vert \nabla u^{\alpha }\left(
\cdot ,\tau \right) \right\Vert _{L^{\infty }}\right) d\tau $. Hence by %
\eqref{eq:grad_v_L_inf_bnd}%
\begin{align*}
\left\Vert w^{\alpha }\left( \cdot ,t\right) \right\Vert _{L^{2}}& \leq
C\alpha ^{3/2}\left\Vert q^{in}\right\Vert _{\dot{B}_{2,\infty
}^{1/2}}\left\Vert \nabla v^{in}\right\Vert _{L^{\infty
}}e^{Ct}e^{C\left\Vert \nabla v^{in}\right\Vert _{L^{\infty }}e^{Ct}} \\
& \leq \alpha ^{3/2}Ce^{Ce^{Ct}},
\end{align*}%
where $C$ depends only on the initial data $q^{in}$.
\end{proof}

\section*{Appendix}

\setcounter{section}{1} \setcounter{theorem}{0} \renewcommand{\thesection}{%
\Alph{section}}

\subsection{\label{sec:GlobalReg}Global regularity of Contour Dynamics-$%
\protect\alpha $ equation}

We consider the vortex patch problem, i.e., a system in which the initial
vorticity $q^{in}$ is proportional to the characteristic function of a
bounded domain $\Omega ^{in}$, $q^{in}=q_{0}\chi _{\Omega ^{in}}$, under the
evolution of the Euler-$\alpha $ equations \eqref{grp:EulerAlphaEqVortForm}.
Due to the conservation of the two-dimensional Euler-$\alpha $ vorticity
along particle trajectories for such an initial data (in fact, it is enough
for the initial vorticity to be in the space of Radon measures in $\mathbb{R}%
^{2}$, see \cite{a_OS01}), the vorticity $q^{\alpha }\left( t\right) $
remains a characteristic function of an evolving in time domain $\Omega
^{\alpha }\left( t\right) $. In this section we present the result which
states that the boundary of a vortex patch evolving under Euler-$\alpha $
equations \eqref{grp:EulerAlphaEqVortForm} remains as smooth, for all time,
as initially boundary, provided the latter is smooth enough in a sense
specified in Theorem \ref{thm:GlobalEx} below.

In two dimensions the evolution of the boundary $x\left( \sigma ,t\right) $
of a vortex patch under the Euler-$\alpha $ equation %
\eqref{grp:EulerAlphaEqVortForm} is given by

\begin{align}
& \frac{\partial x}{\partial t}\left( \sigma ,t\right)
=-q_{0}\int_{S^{1}}\Psi ^{\alpha }\left( \left\vert x\left( \sigma ,t\right)
-x\left( \sigma ^{\prime },t\right) \right\vert \right) \frac{\partial x}{%
\partial \sigma }\left( \sigma ^{\prime },t\right) d\sigma ^{\prime },
\label{eq:CD_alpha} \\
& x\left( \sigma ,0\right) =x^{in}\left( \sigma \right) ,  \notag
\end{align}%
where%
\begin{equation*}
\Psi ^{\alpha }\left( r\right) =\frac{1}{2\pi }\left[ K_{0}\left( \frac{r}{%
\alpha }\right) +\log r\right] ,
\end{equation*}%
and $K_{0}$ denote the modified Bessel functions of the second kind of order
zero. As we mentioned before, for details on Bessel functions, see, e.g.,
\cite{b_W44}. The integro-differential equation \eqref{eq:CD_alpha} is an
analogue of the so-called contour dynamics (CD) equation \cite%
{a_ZHR79,b_MB02}, describing the evolution of the boundary of vortex patch
evolving by means of the Euler equations \eqref{grp:EulerEqVortForm}:
\begin{align}
& \frac{\partial x}{\partial t}\left( \sigma ,t\right) =-\frac{q_{0}}{2\pi }%
\int_{S^{1}}\log \left( \left\vert x\left( \sigma ,t\right) -x\left( \sigma
^{\prime },t\right) \right\vert \right) \frac{\partial x}{\partial \sigma }%
\left( \sigma ^{\prime },t\right) d\sigma ^{\prime },  \label{eq:CDE} \\
& x\left( \sigma ,0\right) =x^{in}\left( \sigma \right) .  \notag
\end{align}%
The CD-$\alpha $ equation \eqref{eq:CD_alpha} is derived following arguments
similar to those used in the Euler case, see \cite{a_ZHR79,b_MB02}.

We show that CD-$\alpha $ equation \eqref{eq:CD_alpha} is well-posed in the
space of Lipschitz functions and in the H\"{o}lder space $C^{n,\beta }$, $%
n\geq 1$, which is the space of $n$-times differentiable functions with H%
\"{o}lder continuous $n^{\text{th}}$ derivative. Let us first describe the H%
\"{o}lder space $C^{n,\beta }\left( S^{1}\subset \mathbb{R};\mathbb{R}%
^{2}\right) $,$\beta \in \left( 0,1\right] $, which is the space of
functions $x:S^{1}\subset \mathbb{R}\rightarrow \mathbb{R}^{2}$, with a
finite norm
\begin{equation*}
\left\Vert x\right\Vert _{n,\beta }=\sum_{k=0}^{n}\left\Vert \frac{d^{k}}{d{%
\sigma }^{k}}x\right\Vert _{C^{0}\left( S^{1}\right) }+\left\vert \frac{d^{n}%
}{d{\sigma }^{n}}x\right\vert _{\beta },
\end{equation*}%
where
\begin{equation*}
\left\Vert x\right\Vert _{C^{0}\left( S^{1}\right) }=\sup_{\sigma \in
S^{1}}\left\vert x\left( \sigma \right) \right\vert
\end{equation*}%
and $\left\vert \cdot \right\vert _{\beta }$ is the H\"{o}lder semi-norm
\begin{equation*}
\left\vert x\right\vert _{\beta }=\sup_{\substack{ \sigma ,\sigma ^{\prime
}\in S^{1}  \\ \sigma \neq \sigma ^{\prime }}}\frac{\left\vert x\left(
\sigma \right) -x\left( \sigma ^{\prime }\right) \right\vert }{\left\vert
\sigma -\sigma ^{\prime }\right\vert ^{\beta }}.
\end{equation*}%
The Lipschitz space $\mathrm{Lip}\left( S^{1}\right) $ is the $C^{0,1}$
space, that is, with the finite norm $\left\Vert x\right\Vert _{\mathrm{Lip}%
\left( S^{1}\right) }=\left\Vert x\right\Vert _{C^{0}\left( S^{1}\right)
}+\left\vert x\right\vert _{1}$. We also use the notation
\begin{equation*}
\left\vert x\right\vert _{\ast }=\inf_{_{\substack{ \sigma ,\sigma ^{\prime
}\in S^{1}  \\ \sigma \neq \sigma ^{\prime }}}}\frac{\left\vert x\left(
\sigma \right) -x\left( \sigma ^{\prime }\right) \right\vert }{\left\vert
\sigma -\sigma ^{\prime }\right\vert }.
\end{equation*}

We consider the CD-$\alpha $ equation \eqref{eq:CD_alpha} as an evolution
functional equation on either one of the following Banach spaces: $\mathrm{%
Lip}$, $C^{1,\beta }$, $\beta \in \left[ 0,1\right] $, $C^{2,\gamma }$, $%
\gamma \in \left[ 0,1\right) $, or $C^{n,\theta }$, $n\geq 3$, $\theta \in
\left( 0,1\right) $. We have the following result

\begin{theorem}
\label{thm:GlobalEx} Let $V$ be either one of the following spaces: $\mathrm{%
Lip}\left( S^{1}\right) $, or $C^{1,\beta }\left( S^{1}\right) $, $\beta \in %
\left[ 0,1\right] $, or $C^{2,\gamma }\left( S^{1}\right) $, $\gamma \in %
\left[ 0,1\right) $, or $C^{n,\theta }\left( S^{1}\right) $, $n\geq 3$, $%
\theta \in \left( 0,1\right) $. Let $x^{in}\in V\cap \left\{ \left\vert
x\right\vert _{\ast }>0\right\} $, then there exists a unique solution {$%
x\in C^{1}\left( \left( -\infty ,\infty \right) ;V\cap \left\{ \left\vert
x\right\vert _{\ast }>0\right\} \right) $} of \eqref{eq:CD_alpha} with
initial value $x\left( \sigma ,0\right) =x^{in}\left( \sigma \right) $. In
particular, if $x_{0}\in C^{\infty }\left( S^{1}\right) \cap \left\{
\left\vert x\right\vert _{\ast }>0\right\} $ then $x\in C^{1}\left( \left(
-\infty ,\infty \right) ;C^{\infty }\left( S^{1}\right) \cap \left\{
\left\vert x\right\vert _{\ast }>0\right\} \right) $.
\end{theorem}

We remark that, although the kernel $\Psi ^{\alpha }$ and its first
derivative $D\Psi ^{\alpha }$ are continuous bounded functions, its higher
derivatives $D^{m}\Psi ^{\alpha }$, $m\geq 2$, are unbounded near the
origin, and the chord arc condition $|x|_{\ast }>0$, which implies simple
curves, allows us to show the integrability of the relevant terms.

We only sketch the main steps of the proof, since it is in the spirit of
\cite{a_BLT08,a_BLT09} and \cite[Chapter 8]{b_MB02}, which can be consulted
for details of such a proof. In \cite{a_BLT09} we show the well-posedness of
vortex sheet problem for Euler-$\alpha $ equations, and it contains various
estimates involving the derivatives of the kernel $\Psi ^{\alpha }$, and
\cite[Chapter 8]{b_MB02} describes the proof of the original vortex patch
problem in the Euler equations case.

The following are the main steps involved in the proof of Theorem \ref%
{thm:GlobalEx}. In the first step, we apply the Contraction Mapping
Principle to the CD-$\alpha $ equation \eqref{eq:CD_alpha} to prove the
short time existence and uniqueness of solutions in the appropriate space of
functions. Next, we derive an \textit{a priori} bound for the controlling
quantity for continuing the solution for all time. At step three one can
extend the result for higher derivatives, using the estimates derived in
\cite{a_BLT09} and \cite[Chapter 8]{b_MB02}.

\subsubsection{Local existence of Contour Dynamics-$\protect\alpha $ equation%
}

Next we show the local existence and uniqueness of solutions in the
Lipschitz space, the details and the estimates in other appropriate spaces
can be done in the same spirit following arguments presented in \cite%
{a_BLT09} and \cite[Chapter 8]{b_MB02}. First we recall some properties of
the kernel $\Psi ^{\alpha }$, see also \eqref{eq:DPsi}-%
\eqref{eq:Der_psi_at_origin}. For $\frac{r}{\alpha }\rightarrow 0$
\begin{equation}
\Psi ^{\alpha }\left( r\right) =\frac{1}{2\pi }\log \alpha +O\left( 1\right)
,  \label{eq:Psi_alpha_near_origin}
\end{equation}%
and for large $\frac{r}{\alpha }$%
\begin{equation}
\Psi ^{\alpha }\left( r\right) =O\left( \log r\right) .
\label{eq:Psi_alpha_at_inf}
\end{equation}%
Also, $D\Psi ^{\alpha }$ is bounded for all $r\in \left[ 0,\infty \right) $,%
\begin{equation}
D\Psi ^{\alpha }\left( r\right) =O\left( \frac{1}{\alpha }\right) ,
\label{eq:Dpsi_bound}
\end{equation}%
where in the big $O$ the constants are independent of $\alpha $. To apply
the Contraction Mapping Principle to the CD-$\alpha $ equation %
\eqref{eq:CD_alpha} we first prove the following result:

\begin{proposition}
\label{prop:cd_alpha_Lip_map}Let $1<M<\infty $, and let $K^{M}$ be the set
\begin{equation*}
K^{M}=\left\{ x\in \mathrm{Lip}\left( S^{1}\right) :\left\Vert x\right\Vert
_{\mathrm{Lip}}<M,\left\vert x\right\vert _{\ast }>\frac{1}{M}\right\} .
\end{equation*}%
Then the mapping
\begin{equation}
x\left( \Gamma \right) \mapsto u\left( x\left( \Gamma \right) \right)
=\int_{S^{1}}K^{\alpha }\left( x\left( \Gamma \right) -x\left( \Gamma
^{\prime }\right) \right) d\Gamma ^{\prime }  \label{eq:vel_map}
\end{equation}%
defines a locally Lipschitz continuous map from $K^{M}$, equipped with the
topology induced by the $\left\Vert \cdot \right\Vert _{\mathrm{Lip}}$ norm,
into $\mathrm{Lip}$.
\end{proposition}

\begin{proof}
We start by showing that $u\left( x\left( \sigma \right) \right) $ maps $%
K^{M}$ into $\mathrm{Lip}$. Let $x\in K^{M}$. By %
\eqref{eq:Psi_alpha_near_origin} and \eqref{eq:Psi_alpha_at_inf} we have%
\begin{align}
\left\vert u\left( x\left( \sigma \right) \right) \right\vert & \leq
q_{0}\int_{S^{1}}\Psi ^{\alpha }\left( \left\vert x\left( \sigma \right)
-x\left( \sigma ^{\prime }\right) \right\vert \right) \left\vert \frac{dx}{%
d\sigma }\left( \sigma ^{\prime }\right) \right\vert d\sigma ^{\prime }
\label{eq:Ux_C0_bnd} \\
& \leq Cq_{0}\left( \log \alpha +1+\log M\right) \left\vert x\right\vert
_{1}.  \notag
\end{align}%
To show Lipschitz continuity of $u\left( x\left( \sigma \right) \right) $ we
use that by mean value theorem and \eqref{eq:Dpsi_bound}, we have that for $%
x\left( \sigma ^{\prime \prime }\right) \in B\left( x\left( \sigma \right)
,\left\vert x\left( \sigma \right) -x\left( \bar{\sigma}\right) \right\vert
\right) $, the ball centered at $x\left( \sigma \right) $ with the radius $%
\left\vert x\left( \sigma \right) -x\left( \bar{\sigma}\right)
\right\vert $,
\begin{align}
\left\vert \Psi ^{\alpha }\left( \left\vert x\left( \sigma \right) -x\left(
\sigma ^{\prime }\right) \right\vert \right) -\Psi ^{\alpha }\left(
\left\vert x\left( \bar{\sigma}\right) -x\left( \sigma ^{\prime }\right)
\right\vert \right) \right\vert & \leq D\Psi ^{\alpha }\left( \left\vert
x\left( \sigma ^{\prime \prime }\right) -x\left( \sigma ^{\prime }\right)
\right\vert \right) \left\vert x\left( \sigma \right) -x\left( \bar{\sigma}%
\right) \right\vert   \label{eq:Psi_alpha_Lip} \\
& \leq \frac{C}{\alpha }\left\vert x\left( \sigma \right) -x\left( \bar{%
\sigma}\right) \right\vert ,  \notag
\end{align}%
and hence
\begin{align*}
\left\vert u\left( x\left( \sigma \right) \right) -u\left( x\left( \bar{%
\sigma}\right) \right) \right\vert & \leq q_{0}\int_{S^{1}}\left\vert \Psi
^{\alpha }\left( \left\vert x\left( \sigma \right) -x\left( \sigma ^{\prime
}\right) \right\vert \right) -\Psi ^{\alpha }\left( \left\vert x\left( \bar{%
\sigma}\right) -x\left( \sigma ^{\prime }\right) \right\vert \right)
\right\vert \left\vert \frac{dx}{d\sigma }\left( \sigma ^{\prime }\right)
\right\vert d\sigma ^{\prime } \\
& \leq \frac{C}{\alpha }q_{0}\left\vert x\right\vert _{1}\left\vert x\left(
\sigma \right) -x\left( \bar{\sigma}\right) \right\vert .
\end{align*}%
Now, we show that $u\left( x\right) $ is locally Lipschitz continuous on $%
K^{M}$. It is enough to prove that for $x\in K^{M}$, $y\in \mathrm{Lip}%
\left( S^{1}\right) $%
\begin{equation}
\left\Vert D_{x}u\left( x\right) y\right\Vert _{\mathrm{Lip}}\leq C\left(
\frac{1}{\alpha },M,\left\Vert x\right\Vert _{\mathrm{Lip}}\right)
\left\Vert y\right\Vert _{\mathrm{Lip}}.  \label{eq:DxUy_Lip_bnd}
\end{equation}%
Let $x\in K^{M}$, $y\in H^{1}\left( S^{1}\right) $, we compute
\begin{align*}
D_{x}u\left( x\left( \sigma \right) \right) y\left( \sigma \right) & =\left.
\frac{d}{d\varepsilon }u\left( x\left( \sigma \right) +\varepsilon y\left(
\sigma \right) \right) \right\vert _{\varepsilon =0} \\
& =-q_{0}\int_{S^{1}}D\Psi ^{\alpha }\left( \left\vert x\left( \sigma
\right) -x\left( \sigma ^{\prime }\right) \right\vert \right) \frac{\left(
x\left( \sigma \right) -x\left( \sigma ^{\prime }\right) \right) \cdot
\left( y\left( \sigma \right) -y\left( \sigma ^{\prime }\right) \right) }{%
\left\vert x\left( \sigma \right) -x\left( \sigma ^{\prime }\right)
\right\vert }\frac{dx}{d\sigma }\left( \sigma ^{\prime }\right) d\sigma
^{\prime } \\
& -q_{0}\int_{S^{1}}\Psi ^{\alpha }\left( \left\vert x\left( \sigma \right)
-x\left( \sigma ^{\prime }\right) \right\vert \right) \frac{dy}{d\sigma }%
\left( \sigma ^{\prime }\right) d\sigma ^{\prime } \\
& =F_{1}\left( x\left( \sigma \right) \right) y\left( \sigma \right)
+F_{2}\left( x\left( \sigma \right) \right) y\left( \sigma \right) .
\end{align*}%
Next we show \eqref{eq:DxUy_Lip_bnd}. To estimate the $C^{0}$ norm we use %
\eqref{eq:Dpsi_bound} for $F_{1}$
\begin{align*}
\left\vert F_{1}\left( x\left( \sigma \right) \right) y\left( \sigma \right)
\right\vert & \leq q_{0}\int_{S^{1}}D\Psi ^{\alpha }\left( \left\vert
x\left( \sigma \right) -x\left( \sigma ^{\prime }\right) \right\vert \right)
\left\vert y\left( \sigma \right) -y\left( \sigma ^{\prime }\right)
\right\vert \left\vert \frac{dx}{d\sigma }\left( \sigma ^{\prime }\right)
\right\vert d\sigma ^{\prime } \\
& \leq \frac{C}{\alpha }q_{0}\left\Vert y\right\Vert _{C^{0}}\left\vert
x\right\vert _{1}.
\end{align*}%
and \eqref{eq:Psi_alpha_near_origin} and \eqref{eq:Psi_alpha_at_inf} for $%
F_{2}$%
\begin{align*}
\left\vert F_{2}\left( x\left( \sigma \right) \right) y\left( \sigma \right)
\right\vert & \leq q_{0}\int_{S^{1}}\left\vert \Psi ^{\alpha }\left(
\left\vert x\left( \sigma \right) -x\left( \sigma ^{\prime }\right)
\right\vert \right) \right\vert \left\vert \frac{dy}{d\sigma }\left( \sigma
^{\prime }\right) \right\vert d\sigma ^{\prime } \\
& \leq Cq_{0}\left( \log \alpha +1+\log M\right) \left\vert y\right\vert
_{1}.
\end{align*}%
Next we show Lipschitz continuity of $D_{x}u\left( x\left( \sigma \right)
\right) y\left( \sigma \right) $. For $F_{2}$ one uses %
\eqref{eq:Psi_alpha_Lip}. For $F_{1}$ we have%
\begin{align*}
|F_{1}\left( x\left( \sigma \right) \right) y\left( \sigma \right) &
-F_{1}\left( x\left( \bar{\sigma}\right) \right) y\left( \bar{\sigma}\right)
|\leq  \\
& \leq q_{0}\int_{S^{1}}\left\vert D\Psi ^{\alpha }\left( \left\vert x\left(
\sigma \right) -x\left( \sigma ^{\prime }\right) \right\vert \right) -D\Psi
^{\alpha }\left( \left\vert x\left( \bar{\sigma}\right) -x\left( \sigma
^{\prime }\right) \right\vert \right) \right\vert \left\vert y\left( \sigma
\right) -y\left( \sigma ^{\prime }\right) \right\vert \left\vert \frac{dx}{%
d\sigma }\left( \sigma ^{\prime }\right) \right\vert d\sigma ^{\prime } \\
& +2q_{0}\int_{S^{1}}D\Psi ^{\alpha }\left( \left\vert x\left( \bar{\sigma}%
\right) -x\left( \sigma ^{\prime }\right) \right\vert \right) \left\vert
x\left( \bar{\sigma}\right) -x\left( \sigma \right) \right\vert \frac{%
\left\vert y\left( \sigma \right) -y\left( \sigma ^{\prime }\right)
\right\vert }{\left\vert x\left( \bar{\sigma}\right) -x\left( \sigma
^{\prime }\right) \right\vert }\left\vert \frac{dx}{d\sigma }\left( \sigma
^{\prime }\right) \right\vert d\sigma ^{\prime } \\
& +q_{0}\int_{S^{1}}D\Psi ^{\alpha }\left( \left\vert x\left( \bar{\sigma}%
\right) -x\left( \sigma ^{\prime }\right) \right\vert \right) \left\vert
y\left( \sigma \right) -y\left( \bar{\sigma}\right) \right\vert \left\vert
\frac{dx}{d\sigma }\left( \sigma ^{\prime }\right) \right\vert d\sigma
^{\prime } \\
& =I_{1}+I_{2}+I_{3}.
\end{align*}%
For $I_{1}$, by the mean value theorem, \eqref{eq:Der_psi_at_origin} and due
to the fact that $\left\vert x\right\vert _{\ast }>\frac{1}{M}$, we have
that for $\sigma ^{\prime \prime }\in S^{1}$ and such that $x\left( \sigma
^{\prime \prime }\right) \in B\left( x\left( \sigma \right) ,\left\vert
x\left( \sigma \right) -x\left( \bar{\sigma}\right) \right\vert \right) $%
\begin{align*}
\left\vert D\Psi ^{\alpha }\left( \left\vert x\left( \sigma \right) -x\left(
\sigma ^{\prime }\right) \right\vert \right) -D\Psi ^{\alpha }\left(
\left\vert x\left( \bar{\sigma}\right) -x\left( \sigma ^{\prime }\right)
\right\vert \right) \right\vert & \leq \left\vert D^{2}\Psi ^{\alpha }\left(
\left\vert x\left( \sigma ^{\prime \prime }\right) -x\left( \sigma ^{\prime
}\right) \right\vert \right) \right\vert \left\vert x\left( \sigma \right)
-x\left( \bar{\sigma}\right) \right\vert  \\
& \leq \left( \frac{1}{4\pi }\frac{1}{\alpha ^{2}}\left\vert \log \frac{%
\left\vert x\left( \sigma ^{\prime \prime }\right) -x\left( \sigma ^{\prime
}\right) \right\vert }{\alpha }\right\vert +\frac{C}{\alpha ^{2}}\right)
\left\vert x\left( \sigma \right) -x\left( \bar{\sigma}\right) \right\vert
\\
& \leq C\left( M\right) \frac{1}{\alpha ^{2}}\left\vert \sigma -\bar{\sigma}%
\right\vert \left( \left\vert \log \left( \frac{\left\vert \sigma ^{\prime
\prime }-\sigma ^{\prime }\right\vert }{\alpha }\right) \right\vert
+1\right) .
\end{align*}%
Therefore,
\begin{align*}
I_{1}& \leq \left\vert \sigma -\bar{\sigma}\right\vert C\left( M\right)
\frac{1}{\alpha ^{2}}q_{0}\left\Vert y\right\Vert _{C^{0}}\left\vert
x\right\vert _{1}\int_{S^{1}}\left( \left\vert \log \left( \frac{\left\vert
\sigma ^{\prime \prime }-\sigma ^{\prime }\right\vert }{\alpha }\right)
\right\vert +1\right) d\sigma ^{\prime } \\
& \leq C\left( M,\frac{1}{\alpha },q_{0}\right) \left\Vert y\right\Vert
_{C^{0}}\left\vert \sigma -\bar{\sigma}\right\vert .
\end{align*}%
For $I_{2}$ and $I_{3}$ we use \eqref{eq:Dpsi_bound} and $\left\vert
x\right\vert _{\ast }>\frac{1}{M}$ to obtain
\begin{equation*}
I_{2},I_{3}\leq C\left( M,\frac{1}{\alpha },q_{0}\right) \left\vert
y\right\vert _{1}\left\vert \sigma -\bar{\sigma}\right\vert .
\end{equation*}
\end{proof}

Proposition \ref{prop:cd_alpha_Lip_map} implies the local existence and
uniqueness of solutions:

\begin{proposition}
Let $K^{M}=\left\{ x\in \mathrm{Lip}\left( S^{1}\right) :\left\Vert
x\right\Vert _{\mathrm{Lip}}<M,\left\vert x\right\vert _{\ast }>\frac{1}{M}%
\right\} $and let $x_{0}\in \mathrm{Lip}\left( S^{1}\right) \cap \left\{
\left\vert x\right\vert _{\ast }>0\right\} $, then for any $M$, $1<M<\infty $%
, such that $x_{0}\in K^{M}$, there exists a time $T(M)$, such that the
system \eqref{eq:CD_alpha} has a unique local solution $x\in
C^{1}((-T(M),T(M));K^{M})$.
\end{proposition}

\subsubsection{Global existence of Contour Dynamics-$\protect\alpha $
equation}

To show the global existence of the CD-$\alpha $ equation, we assume by
contradiction, that $T_{\max }<\infty $, where $\left[ 0,T_{\max }\right) $
is the maximal interval of existence, and hence the solution leaves in a
finite time the open set $K^{M}$, for all $M>1$, that is, $%
\limsup_{t\rightarrow T_{\max }^{-}}\left\Vert x\right\Vert _{V}=\infty $ or
$\limsup_{t\rightarrow T_{\max }^{-}}\frac{1}{\left\vert x\left( \cdot
,t\right) \right\vert _{\ast }}=\infty $. Therefore, if we show global
bounds on $\frac{1}{\left\vert x\left( \cdot ,t\right) \right\vert _{\ast }}$
and $\left\Vert x\left( \cdot ,t\right) \right\Vert _{V}$ in $\left[
0,T_{\max }\right) $, we obtain a contradiction to the blow-up, and thus the
obtained local solutions can be continued for all time. The result extends
to negative times as well.

To control the quantities $\frac{1}{\left\vert x\left( \cdot ,t\right)
\right\vert _{\ast }}$ and $\left\Vert x\left( \cdot ,t\right) \right\Vert
_{V}$ one needs to bound $\int_{0}^{T_{\max }}\left\Vert \nabla
_{x}u^{\alpha }\left( x(\cdot ,t),t\right) \right\Vert _{L^{\infty }}dt$.
Next proposition shows the bound on $\left\Vert \nabla u^{\alpha
}\right\Vert _{L^{\infty }}$ for the vortex patch initial data.

\begin{proposition}
Let $q^{in}\in L^{1}\left( \mathbb{R}^{2}\right) \cap L^{\infty }\left(
\mathbb{R}^{2}\right) $ then%
\begin{equation*}
\left\Vert \nabla _{x}u^{\alpha }\left( x(\cdot ,t),t\right) \right\Vert
_{L^{\infty }}\leq C\left( \frac{1}{\alpha }\right) \left( \left\Vert
q^{in}\right\Vert _{L^{\infty }}+\left\Vert q^{in}\right\Vert
_{L^{1}}\right) .
\end{equation*}
\end{proposition}

\begin{proof}
We write
\begin{align*}
\nabla _{x}{u}^{\alpha }\left( x,t\right) & =\int_{\mathbb{R}^{2}}\nabla
K^{\alpha }\left( x-y\right) q^{\alpha }\left( y,t\right) dy \\
& =\int_{\left\vert x-y\right\vert <\alpha }+\int_{\left\vert x-y\right\vert
\geq\alpha }=I_{1}+I_{2}.
\end{align*}%
Using \eqref{eq:Der_psi_at_origin} we obtain
\begin{align*}
I_{1}& \leq C\left\Vert q^{\alpha }\left( \cdot ,t\right) \right\Vert
_{L^{\infty }}\int_{0}^{\alpha }\left\vert D^{2}\Psi ^{\alpha }\left(
r\right) \right\vert rdr \\
& \leq C\left\Vert q^{in}\right\Vert _{L^{\infty }}\int_{0}^{\alpha
}r\left\vert \frac{1}{4\pi }\frac{1}{\alpha ^{2}}\log \frac{r}{\alpha }+%
\frac{C}{\alpha ^{2}}\right\vert dr \\
& \leq C \left\Vert q^{in}\right\Vert _{L^{\infty }}
\end{align*}%
and%
\begin{align*}
I_{2}& \leq \sup_{\left\vert x-y\right\vert \geq\alpha }\left\vert \nabla
K^{\alpha }\left( x-y\right) \right\vert \int_{\mathbb{R}^{2}}q^{\alpha
}\left( y,t\right) dy \\
& \leq C\frac{1}{\alpha ^{2}}\left\Vert q^{in}\right\Vert _{L^{1}}.
\end{align*}
\end{proof}

\section*{Acknowledgements}

This work was supported in part by the BSF grant no.~2004271, the
ISF grant no.~120/06, and the NSF grant no.~DMS-0708832.

\bibliographystyle{siam}
\bibliography{VortexSheetBib}

\end{document}